  \documentclass[final,3p,times]{elsarticle}

\usepackage{amssymb}

\usepackage{amsthm}

\theoremstyle{definition}

\usepackage{amsmath}

\usepackage{lineno}
\usepackage{algorithm}
\usepackage[noend]{algpseudocode}
\usepackage{booktabs} 
\usepackage{color}
\usepackage{tikz}
\usetikzlibrary{calc}
\usetikzlibrary{positioning}
\usepackage{placeins} 
\usepackage{subcaption}
\usepackage{graphicx}

\usepackage{threeparttable}

\newcommand{\myset}[3]{\{#1_{#2}\}^{#3}_{#2=1}}
\newcommand{\Real}[1]{\operatorname{Re}\left\{ #1 \right\}}
\newcommand{\Imag}[1]{\operatorname{Im}\left\{ #1 \right\}}

\makeatletter
\def\BState{\State\hskip-\ALG@thistlm}
\makeatother

\makeatletter
\def\ps@pprintTitle{%
  \let\@oddhead\@empty
  \let\@evenhead\@empty
  \def\@oddfoot{\reset@font\hfil\thepage\hfil}
  \let\@evenfoot\@oddfoot
}
\makeatother

\begin{document}

\begin{frontmatter}



\title{Partition of Unity Extension of Functions on Complex Domains}


\author[KTH]{Fredrik Fryklund\corref{cor1}}
\ead{ffry@kth.se}
\cortext[cor1]{Corresponding author}
\author[KTH]{Erik Lehto}
\ead{elehto@kth.se}
\author[KTH]{Anna--Karin Tornberg}
\ead{akto@kth.se}

\address[KTH]{KTH Mathematics, Swedish e--Science Research Centre,\\ 100 44 Stockholm, Sweden}
\begin{abstract}

We introduce an efficient algorithm, called partition of unity
extension or PUX, to construct an extension of desired regularity of a
function given on a complex multiply connected domain in $2D$.
Function extension plays a fundamental role in extending the
applicability of boundary integral methods to inhomogeneous partial
differential equations with embedded domain techniques. Overlapping
partitions are placed along the boundaries, and a local extension of
the function is computed on each patch using smooth radial basis
functions; a trivially parallel process. A partition of unity method
blends the local extrapolations into a global one, where weight
functions impose compact support.  The regularity of the extended function can be controlled by the
construction of the partition of unity function.  We evaluate the
performance of the PUX method in the context of solving the Poisson
equation on multiply connected domains using a boundary integral
method and a spectral solver. With a suitable choice of parameters the error converges as a tenth order method down to $10^{-14}$.

\end{abstract}


\end{frontmatter}

 
\section{Introduction}
\label{S:1}
This paper addresses the issue of how to numerically construct an
extension of a function defined on a complex domain. Without prior acquaintance with the topic it could appear simple.
However, adding requirements on the global regularity of the extended
function it becomes a non-trivial task. Furthermore, it is often
desirable that the extended function has compact support and that it
can be efficiently constructed.   One important application for function extension is to extend
the applicability of integral equation methods for solving partial differential equations (PDEs).
Integral equation methods have been shown to be both highly accurate
and efficient when solving homogeneous constant coefficient elliptic
partial differential equations in complex geometry. Function extension
is a key component in a framework for solving non-homogeneous elliptic
PDEs, and furthermore to solve time-dependent equations such as the heat
equation and extending from the solution of Stokes equations to
Navier-Stokes equations \citep{Kropinski2011modHelm,Kropinski2011Heat,doi:10.1137/S1064827597317648}. \\

The idea is to avoid solving the full linear elliptic inhomogeneous PDE, given on a complex domain, by
splitting the problem into two. The right–-hand side is extended to a geometrically
simpler domain, such as a box, and as part of the full solution, a
particular solution is computed on this simpler domain. The
homongeneous problem is solved on the original domain with modified
boundary conditions, such that the total solution is the sum of the
two. This general idea  has been used also for other numerical
methods and this group of methods is often referred  to as embedded
boundary techniques. For simple geometries an arsenal of powerful solution methods are
available, but their accuracy is often limited by the global
regularity of the extended function. Several different function
extension methods, used in this context, have been suggested in recent
years \cite{ASKHAM20171,Albin2011,Stein2017155,Bruno20102009,Lyon20103358,Shirokoff2015}.\\

The approach to function extension suggested by Askham et al. \cite{ASKHAM20171} and
Stein et al. \cite{Stein2017155} are both based on the same idea, but quite different in implementation. Given $f$ on a domain
$\bar{\Omega}$, use the vaules of $f$ on the boundary of $\Omega$ as the Dirichlet
data for the external Laplace problem
  \begin{alignat}{2}
\label{eq:laplace_askham}
\Delta w & {}={}   0 {} \text{ in }  \mathbb{R}^{2}\backslash\ \Omega,  \\
\label{eq:laplace_askhambc}
w & {}={}   f \text{ on }  \partial \Omega.
\end{alignat}
Then a globally continuous extension $f^{e}$ of $f$ is given by 
\begin{equation}
  \label{eq:functionextension_askham}
      f^{e}(\mathbf{y}) = \begin{cases}f(\mathbf{y}),\quad \mathbf{y}\in\Omega,\\
      w(\mathbf{y}),\quad \mathbf{y}\in\mathbb{R}^{2}\backslash\ \Omega.
      \end{cases}
\end{equation} 
In \cite{ASKHAM20171}, this problem is solved by an integral equation based method,
whereas in \cite{Stein2017155} it is coupled with a so called \emph{immersed boundary system of equations}.  This extension will be  in $C^{0}$, but will not have
compact support.  To obtain higher regularity the biharmonic equation
can be solved instead of \eqref{eq:laplace_askham}--\eqref{eq:laplace_askhambc}, or even polyharmonic. Their respective
work show the complexity of the problem and what price is considered
reasonable to pay to obtain an extension. Moreover both observe
that the accuracy of the solution to the associated PDE relies heavily
on the regularity of said extension $f^{e}$.\\

Another alternative, given by \citep{Shirokoff2015}, is to extend the  solution to the PDE, instead of the right hand side, and use an
active penalty method. A function extension of global regularity $k$ is
created by matching normal derivatives of degree $k$ of the given boundary data. The extension is expressed in a basis that is rapidly
decaying with the distance to the boundary. Function extension can
also be achieved by Fourier continuation methods: in $1D$ the domain of
interest is embedded into a larger one and a smooth periodic extension
is constructed, which yields an appropriate setting for spectral
methods. Dimensional splitting is used for higher dimensional problems. See \citep{Lui2009,Bruno2010,Albin2010,ALBIN20116248} and the references therein. The methods
and the associated references included above is by no means a
complete list of methods for function extension. In all mentioned
cases above, no higher than a fourth order method is obtained for solving the Poisson equation.\\

In this paper, we present a new method, \emph{Partition of Unity Extension},
or PUX, to compute a compactly supported extension of a function.  We
assume that the values of a function $f$ are known at all points of a
regular grid that fall inside a domain $\Omega$, and we want to
compute the values of the extended function on this regular grid
outside of $\Omega$. The domain  $\Omega$ can be multiply connected. In the PUX method, overlapping circular partitions or patches are
placed along the boundaries such that each is intersected by the boundary $\partial\Omega$ and
a local extension is defined on each patch. A second layer of patches
is placed outside of the first, on which the local values are defined
to be zero.  These zero patches enter the definition of the partition
of unity function that is used to blend the local extensions into a
global one, imposing compact support and regular decay to zero. The
choice of functions used to build up the partition of unity function
determines the regularity of the extended function.\\

The local extensions on the patches intersected by $\partial\Omega$ are
determined using radial basis functions (RBFs). RBF centres are placed
irregularly with the same distribution for each circular patch, and an
RBF interpolant is determined via a least squares problem, using the
values of $f$ on the regular points inside $\Omega$. The values of the
local extension are then computed on the regular points inside the
patch that fall outside of $\Omega$. Always centring the patches at 
grid points of the regular grid, a matrix $A$ can be precomputed once
and be used for all patches. For each patch, an identification is made
of which points are inside and outside $\Omega$ and a local least
squares problem with the relevant rows of $A$ is solved with the
inside data, a trivially parallel task.\\

To assess the quality of a function extension it must be considered in
its context of use, as there is no unique extension over the boundary
of a domain. In this paper, we use it to solve the Poisson equation,
using an integral equation approach. Thus the results by Askham et al. in \cite{ASKHAM20171} are suitable for comparison. \\

The paper is organised as follows: in section \ref{s:Pequation} we detail how we
solve the Poisson equation assuming an extension of the right hand
side $f$ is known. In section \ref{s:background} we introduce the concepts and
techniques from RBF interpolation and the partition of unity method
that we need to introduce our method. The PUX method is presented in
section \ref{s:functionextesion}, where a function extension is constructed. Thereafter
follows section \ref{s:errorAnalysis} with a discussion of sources of errors associated
with function extension and solving the Poisson equation. Section \ref{s:summary}
is a summary, combined with implementation details, for solving the
Poisson equation with the techniques described in this paper. In
section \ref{s:numericalResults} we perform numerical experiments and carefully discuss
parameter choices.  Finally our conclusions and an outlook are
presented in section \ref{s:conclusions}.

\section{The Poisson equation on irregular domains in a boundary integral method environment}
\label{s:Pequation}
\label{S:2}

To understand why a function extension is useful for solving linear elliptic PDEs, and why its construction is motivated to pursue, we sketch the solution procedure. Consider the Poisson equation with Dirichlet boundary conditions, stated as
\begin{alignat}{4}
\label{eq:poisson}
\Delta u & {}={} & f {} &\text{ in }  \Omega, \\
\label{eq:poissonBC}
u & =  & g &\text{ on } \partial \Omega,
\end{alignat}
where $\Omega$ is a  simply or multiply connected compact domain in $\mathbb{R}^{2}$ with a Lipchitz continuous 
boundary, meaning corners but not cusps are allowed. We will refer to the Poisson equation as the full problem. Introduce homogeneous and 
particular solutions $u^{H}$ and $u^{P}$ such that $u = u^{H}+u^{P}$.  Let $B$ be a simple box domain of size $[-L,L]^{2}$ which embeds $\Omega$, and let $E$ be the complement of $\bar{\Omega}$ relative $B$. We refer to $E$ as the \emph{extension domain}. If we know an extension $f^{e}$ of $f$ from $\Omega$ to $B$,  such that $f^{e} =  f \text{ in }  \Omega$ and $\text{supp}(f^{e})\subset B$, the particular solution $u^{P}$ can be computed from
\begin{alignat}{2}
\label{eq:poissonR2}
\Delta u^{P} & {}={}  f^{e} {} \text{ in }  \mathbb{R}^{2}, \\
\label{eq:poissonR2L2}
u^{P} &\in  L^{2}(\mathbb{R}^{2}) {},
\end{alignat}
assuming $f^{e}$ can be constructed numerically, which is the focus of this paper. The homogeneous solution $u^{H}$ is obtained by solving the Laplace equation with modified Dirichlet boundary conditions, i.e. 
\begin{alignat}{2}
\label{eq:laplace_mod}
\Delta u^{H} & {}={}   0 {} \text{ in }  \Omega,  \\
\label{eq:laplace_modbc}
u^{H} & {}={}   g {}-{} u^{P}\vert_{\partial \Omega} \text{ on }  \partial \Omega.
\end{alignat}
 The numerical treatment of the Laplace equation on $\Omega$  is discussed in subsection \ref{ss:homogeneuous}. It  can be solved for the homogeneous solution to high precision very effectively with a boundary integral method. In our specific setting, the Dirichlet boundary data will be modified according to the solution of the free--space Poisson equation  \eqref{eq:poissonR2}, as given in  \eqref{eq:laplace_modbc}. If the extension $f^{e}$ of $f$ is known and has compact support one can truncate the free--space Poisson equation \eqref{eq:poissonR2}--\eqref{eq:poissonR2L2} to a box; then a numerical solution can efficiently be obtained by fast spectral methods, as described in \ref{ss:particularSolution}. Therefore it is imperative $f^{e}$ has high global regularity. Note that both \eqref{eq:poissonR2}--\eqref{eq:poissonR2L2}  and \eqref{eq:laplace_mod}--\eqref{eq:laplace_modbc} are relatively simple problems to solve numerically to  high accuracy at a low computational cost, if said extension $f^{e}$ is simple to construct. The solution to the full problem \eqref{eq:poisson}--\eqref{eq:poissonBC} is given $u = u^{H} + u^{P}$.\\


\subsection{Homogeneous solution}
\label{ss:homogeneuous}
Consider the Dirichlet problem for the Laplace equation,
\begin{alignat}{1}
\label{eq:laplace}
\Delta u^{H} & =  0 \text{ in }  \Omega,  \\
\label{eq:laplacebc}
u^{H} & =   \tilde{g} \text{ on } \partial \Omega,
\end{alignat}
with $\Omega$ assumed to be simply connected for now. The numerical solution to a linear elliptic homogeneous PDEs is attractive to compute using a boundary integral method, as it can be both accurate and efficient. The boundary integral formulation becomes more dense if expressed with a complex representation, which is natural to adopt since we consider the plane $\mathbb{R}^{2}$. The use of complex notation will be isolated for treating \eqref{eq:laplace}--\eqref{eq:laplacebc} and all expressions have an equivalent counterpart in $\mathbb{R}^{2}$. We will move freely between the two representations as we see fit. In $\mathbb{C}$ points will be denoted $\tau$ and $z$, and the real and imaginary parts of any $\tau\in \mathbb{C}$ are denoted  $\Real{\tau}$ and $\Imag{\tau}$. \\
 
By introducing an unknown density $\mu\!:\!\partial \Omega\rightarrow \mathbb{R}$   the solution $u^{H}$ of \eqref{eq:laplace}--\eqref{eq:laplacebc} can be represented with a double layer potential 

\begin{equation}
    \label{eq:doublelayer}
    u^{H}(z) = \frac{1}{2\pi}\int\limits_{\partial \Omega} \mu(\tau)\Imag{ \frac{d\tau}{\tau -z}  }, \quad \forall z\in\Omega.
\end{equation}

Clearly it is sufficient to know $\mu$ on the boundary $\partial \Omega$ to have an expression for $u^{H}$ anywhere in $\Omega$. Thus only  $\partial \Omega$ needs to be discretised and the dimensionality of \eqref{eq:laplace}--\eqref{eq:laplacebc} has been reduced by one. To obtain the unknown density function $\mu$ we  treat \eqref{eq:doublelayer} in the sense of a Cauchy principal value as we pass the limit $z\rightarrow z_{0}$ for $z \in \Omega$ and apply the boundary
condition. This yields
\begin{equation}
	\label{eq:doublelayerMu}
	\tilde{g}(z_{0}) = \frac{1}{2} \mu(z_{0}) + \frac{1}{2\pi}\int\limits_{\partial \Omega}  \mu(\tau)\Imag{ \frac{d\tau}{\tau -z_{0}}  }, \quad z_{0}\in\partial\Omega.
\end{equation}
The integrand has a well--defined  limit as $\tau\rightarrow z_{0}$. This is a Fredholm integral equation of second kind. If $\partial \Omega$ is Lipschitz continuous the integrand is a 
\emph{compact} integral operator, thus the Fredholm alternative states 
that \eqref{eq:doublelayerMu} is uniquely solvable for $\mu$
\cite{atkinson1997numerical}. This property is inherited by the discretised system, which is obtained by dividing the boundary $\partial\Omega$ into panels of equal length and on each apply the $16$ points \emph{Gauss--Legendre quadrature rule}. The total number of Gauss--Legendre panels  is denoted $N_{\partial \Omega}$. We use  Nystr\"om's method, meaning the collocation points are set to coincide with the  quadrature points $z_{j}$, with $j = 1,\ldots,16\, N_{\partial \Omega}$. The resulting linear system is 

\begin{equation}
\label{eq:doublelayerMu_disc}
     \tilde{g}(z_{i})=\sum\limits_{\substack{j=1\\ j\neq i}}^{16\,N_{\partial\Omega}}\mu(z_{j})\omega_{j}\Imag{ \frac{\tau_{i}'}{\tau_{i} -z_{j}}  }  + \mu(z_{i})\left(\omega_{j}\Imag{ \frac{\tau_{i}''}{\tau_{i}'}  } + \frac{1}{2}\right),\quad i=1,2\ldots,16\,N_{\partial \Omega},
\end{equation}
with $\omega_{j}$ denoting  quadrature weights. Here $\tau_{i}'$ and $\tau_{i}''$ mean that the boundary $\partial \Omega$ has been parametrised as $\tau = \tau(t)$  and differentiated at $t_{i}$ once and twice, respectively.  The corresponding matrix--vector representation of \eqref{eq:doublelayerMu_disc} yields a dense system matrix. A property of the discretisation
of a second kind integral equation is that the condition number does not increase with finer resolution of $\partial \Omega$. We solve for $\mu$ at the quadrature points with \emph{GMRES}. \\

\begin{figure}[ht]
    \begin{center}
        \includegraphics[width=0.3\textwidth]{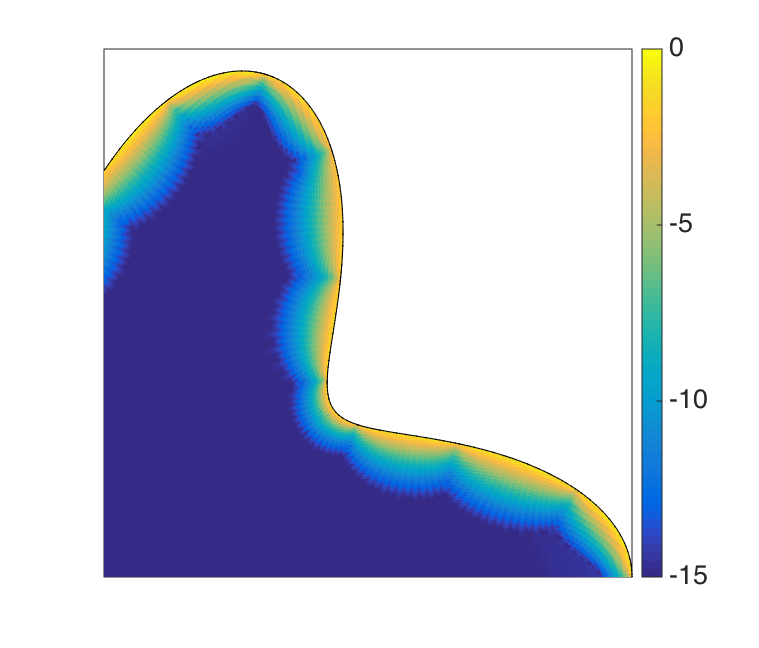}
        \includegraphics[width =0.3\textwidth]{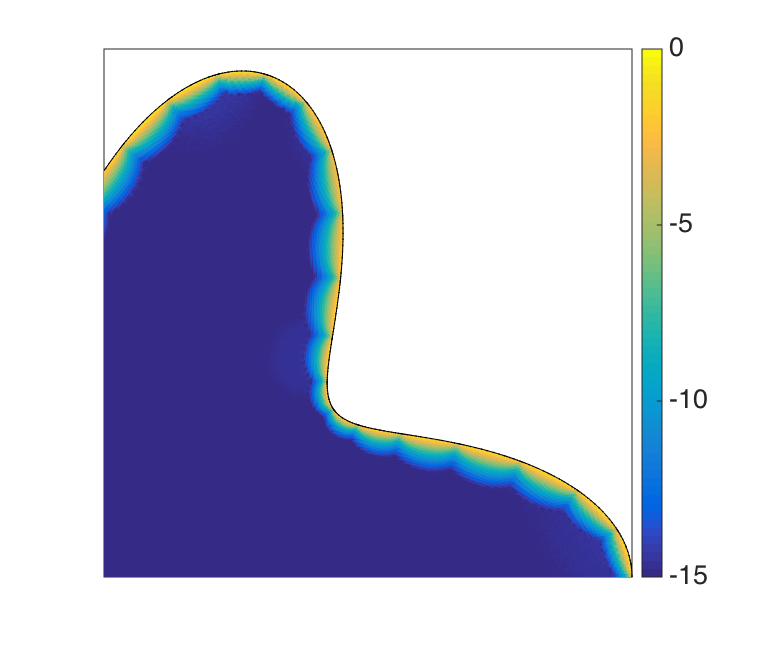}
        \includegraphics[width=0.3\textwidth]{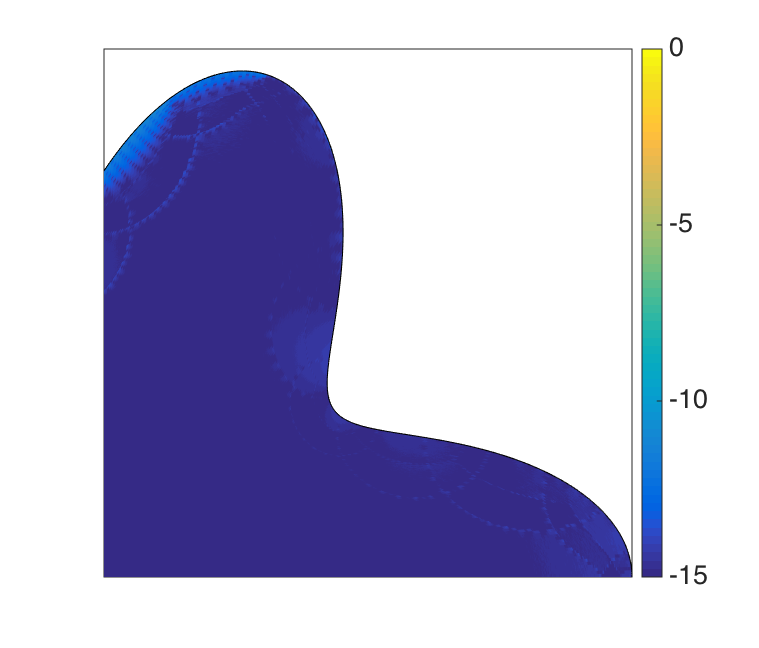}
    \end{center}
    \caption{Example of pointwise relative error, in $\log_{10}$--scale, in part of the computational domain when solving the Laplace equation using:  normal quadrature using $35$ panels (left), normal quadrature using $70$ panels (centre) and special quadrature using $35$ panels (right). }
    \label{fig:normquadspecquad}
\end{figure}

Once $\mu$ is known, $u^{H}$ can be computed at the points $z$ in $\Omega$ of interest by evaluating a numerical approximation of the integral in \eqref{eq:doublelayer}. The integrand, however, becomes nearly singular for $z$ close to $\partial \Omega$, which causes a loss of accuracy if the regular panel based $16$ point Gauss--Legendre rule is used. This can to some extent be remedied by increasing the number of panels along $\partial \Omega$, as in Figure \ref{fig:normquadspecquad} where the centre image has twice as many panels as the left. Since the error increases exponentially fast as the boundary is approached \cite{Klinteberg2017}, this can not provide accurate solutions arbitrarily close to the boundary. Excellent results can however be obtained by applying special quadrature techniques.\\

We use one such technique, namely an interpolatory quadrature method introduced by Ojala and Helsing \cite{specialQuad}. 
This method is based on expanding the complex density $\mu$ as a polynomial in the complex variable $\tau$ over one panel and invoking recursive formulas to analytically evaluate all integrals that are needed. The coefficients of the polynomial 
are the solution of a Vandermonde system. It is however possible to solve the transposed problem instead, for which   the right hand side depends only on the panel geometry and the location of the evaluation point, and not on the discrete values of $\mu$ on the panel. Together with a scaling and rotation of the panel, this controls the ill--conditioning of the problem.
\\

If an evaluation point is within a panel length from a panel's midpoint, we check if the special quadrature method needs to be invoked. Here, we use one of the above mentioned analytically known integrals from expanding $\mu$ and check the accuracy in the numerically obtained value (with the Gauss--Legendre quadrature rule) for this integral. It contains the same near singularity as \eqref{eq:doublelayer} and is a good
indicator of the accuracy that will be obtained by the regular quadrature.\\

So far only simply connected domains have been considered. To fix notation for the Laplace equation on multiply connected domains, let  $\Omega$ consist of a finite $(\kappa +1)$--connected region with a boundary $\partial \Omega$ consisting of $(\kappa +1)$ closed curves. These are denoted $\partial \Omega_{0}$, $\partial \Omega_{1}$, $\ldots$, $\partial \Omega_{\kappa}$, where  $\partial \Omega_{0}$ forms the outer boundary  of the region $\Omega$. Let the \emph{cavity} $\Omega_{k}$ be the set enclosed by $\partial \Omega_{k}$,  disjoint with $\Omega$, for $k = 1,2,\ldots,\kappa$. For such a setting the associated integral equation to \eqref{eq:doublelayer} has a nontrivial nullspace; consequently the Fredholm alternative is not applicable. Thus modifications to the method described above are in order. We apply the remedy suggested in \cite{Greenbaum1993} and the references therein. Instead of \eqref{eq:doublelayer}, the solution to Laplace equation on multiply connected domains is

\begin{equation}
    \label{eq:doublelayerMult}
     u^{H}(z) = \frac{1}{2\pi}\int\limits_{\partial \Omega} \mu(\tau)\Imag{ \frac{d\tau}{\tau -z}  } + \sum\limits_{k=1}^{\kappa}A_{k}\log|z-z_{k}|, \quad \forall z\in\Omega,
\end{equation}
where $z_{k}$ is an interior point of $\Omega_{k}$, and the coefficients $A_{k}$ and the density $\mu$ are unknown. By imposing 

\begin{equation}
    \label{eq:doublelayerMu1Mult}
     \int\limits_{\partial \Omega_{k}} \mu(\tau)\,d|\tau|=0, \quad \text{for } k = 1,2,\ldots,\kappa
\end{equation}
where integration is due to arc length, the resulting boundary integral equation 

\begin{equation}
\begin{cases}
    \label{eq:doublelayerMu2Mult}
    \tilde{g}(z_{0})= \frac{1}{2} \mu(z_{0}) + \frac{1}{2\pi}\int\limits_{\partial \Omega}  \mu(\tau)\Imag{ \frac{d\tau}{\tau -z_{0}}  } +\sum\limits_{k=1}^{\kappa}A_{k}\log|z_{0}-z_{k}|, \quad z_{0}\in\partial\Omega.\\
    \int\limits_{\partial \Omega_{k}}  \mu(\tau)\,d|\tau| =0 , \quad \text{for } k = 1,2,\ldots,\kappa.
    \end{cases}
\end{equation}
satisfies the Fredholm alternative. Hence the coefficients $A_{k}$ and the density $\mu$ can be solved for simultaneously via \eqref{eq:doublelayerMu2Mult}. Observe that the special quadrature techniques for evaluating \eqref{eq:doublelayerMult} are the same as for the simply connected case, as the points $z_{k}$ are chosen so that $\log|z-z_{k}|$ pose no numerical difficulties to evaluate.\\

While not necessary to solve the Laplace equation, a useful aspect of \eqref{eq:doublelayer} is that for $\mu \equiv 1$ one has

\begin{equation}
    \label{eq:doublelayer_mu1}
    \frac{1}{2\pi}\int\limits_{\partial \Omega} \Imag{ \frac{d\tau}{\tau -z}  } = \begin{cases}
    1,& \text{if } z\in\Omega,\\
    1/2,              & \text{if } z\in\partial\Omega,\\
    0,              & \text{if } z\notin\bar{\Omega}.
\end{cases}
\end{equation}

Thus the framework for solving the Laplace equation can be used to identify points as in $\Omega$, on the boundary of $\Omega$ or outside  $\bar{\Omega}$.\\

\subsection{Particular solution}
\label{ss:particularSolution}
The particular solution to the full problem \eqref{eq:poisson}--\eqref{eq:poissonBC} is acquired by solving the free--space Poisson equation \eqref{eq:poissonR2}--\eqref{eq:poissonR2L2}. A key component is an extension $f^{e}$ of $f$; its construction with PUX is described in section \ref{s:functionextesion}. If $f^{e}$ is known, such that $f^{e} = f$ in $\Omega$ and $\text{supp}(f^{e})\subset B$, the solution to \eqref{eq:poissonR2}--\eqref{eq:poissonR2L2} is known to be
\begin{equation}
\label{eq:poissonPartSol}
    u^{P}(\mathbf{x}) = \int\limits_{\mathbb{R}^{2}} K(\|\mathbf{x}-\mathbf{y}\|_{2})f^{e}(\mathbf{y})\,d\mathbf{y} = \int\limits_{B} K(\|\mathbf{x}-\mathbf{y}\|_{2})f^{e}(\mathbf{y})\,d\mathbf{y},
\end{equation}
where $K$ is the Green's function. In $\mathbb{R}^{2}$ one has
\begin{equation}
    K(r) = \frac{-1}{2\pi}\log(r).
\end{equation}
For a derivation see any basic textbook on the subject, e.g. \cite{evans2010partial}.
To handle \eqref{eq:poissonPartSol} numerically we apply the approach suggested by Vico et al. in \cite{MR3535728}: since $f^{e}$ has compact support in $B$, the Green's function can be truncated without changing the value of $u^{P}(\mathbf{x})$ in \eqref{eq:poissonPartSol} for $\mathbf{x}\in B$. Replace $K$ with

\begin{equation}
\label{eq:truncatedGreen}
    \tilde{K}(r) = \begin{cases}
    \frac{-1}{2\pi}\log(r) &\text{for } r < R = \frac{3}{2}L, \\
    0 & \text{ otherwise},
    \end{cases}
\end{equation}
where $[-L,L]^{2}$ is the box $B$ containing $\text{supp}(f)$;  $3/2\,L$ is slightly larger than the greatest possible distance between two points in $B$. Denote the Fourier transform  of $\tilde{K}$ as $\hat{\tilde{K}}$, which has a closed analytical expression, namely 
\begin{equation}
    \hat{\tilde{K}}(k) = \frac{1-J_{0}(Rk)}{k^{2}}-\frac{R\log(R)J_{1}(Rk)}{k},
    \label{eq:truncatedSpectralGreen}
\end{equation}
where $J_{0}$ and $J_{1}$ are the Bessel functions of the first kind of order $0$ and $1$, respectively.  Here $k = \|\mathbf{k}\|_{2}$ with $\mathbf{k}\in \mathbb{R}^{2}$, thus it is radial. 
Furthermore, the limit as $k$ goes to $0$ is well--defined:
\begin{equation}
    \lim\limits_{k \rightarrow 0}\hat{\tilde{K}}(k) = \frac{R^{2}}{4}(1-2\log(R)).
    \label{eq:truncatedSpectralRepofGreen}
\end{equation}
The Fourier coefficients for $u^{P}$ are

\begin{equation}
    \hat{u}^{P}(\mathbf{k}) = \hat{\tilde{K}}(k) \hat{f}^{e}(\mathbf{k}), \quad \mathbf{k}\in\mathbb{R}^{2},
    \label{eq:fouroerCoeff}
\end{equation}
and the solution $u^{P}(\mathbf{x})$ is given by the inverse Fourier transform 
\begin{equation}
\label{eq:freeSol}
u^{P}(\mathbf{x}) = \frac{1}{(2\pi)^{2}}\int\limits_{\mathbb{R}^{2}}\hat{u}^{P}(\mathbf{k})e^{i\mathbf{k}\cdot \mathbf{x}}\,d\mathbf{k}.
\end{equation}
Clearly it is beneficial to have $\hat{\tilde{K}}(k)$ instead of $\hat{K}(k) = -1/k^{2}$, which has a singularity  at $k = 0$. The integral in \eqref{eq:freeSol} can be discretised  with the trapezoidal rule, which allows the usage of an inverse FFT to compute it. It will be sufficiently resolved in $k$--space if upsampling of the FFT is done with a factor greater or equal to $5/2$ for $R = 3/2\,L$ \cite{afKlinteberg2017}. Hence, to numerically obtain $u^{P}$ on a uniform grid in $B$ we apply FFTs to transform the values of $f^{e}$ and $\tilde{K}$ given on a uniform grid. The FFTs are zero padded to obtain the required upsampling, and an inverse FFT is applied after multiplication of the coefficients \eqref{eq:freeSol}. This zero padding avoids pollution from periodic copies of the truncated kernel, and can be reduced to a factor of two after a precomputation step for $\hat{\tilde{K}}$, see \cite{MR3535728}. This is the standard minimum oversampling to compute an aperiodic convolution. This yields an FFT based method for the free--space Poisson equation that converges spectrally as the uniform grid over $B$ is refined, if $f^{e}$ is smooth.  \\

 To solve the full problem \eqref{eq:poisson}--\eqref{eq:poissonBC} the solution $u^{p}$ is needed on $\partial \Omega$ to define the modified boundary conditions for the Laplace equation \eqref{eq:laplace_mod}--\eqref{eq:laplace_modbc}. A non--uniform FFT with the Fourier coefficients $\hat{u}^{P}$obtained above allows us to evaluate $u^{P}$ on $\partial \Omega$ \cite{MR2115056}. \\

\section{Background for interpolation with RBFs and partition of unity methods} 
\label{s:background}
The PUX--method is based on the structure of interpolation  with radial basis functions (RBFs) and partition of unity methods. This section aims to present a comprehensive background for these techniques, before the function extension is discussed in section \ref{s:functionextesion}. An informal description of RBFs and partition of unity methods for the interpolation problem can be summarised as follows. The interpolation domain $\Omega$ is covered with overlapping partitions. On each partition a local interpolant, based on $f$, is constructed with RBFs. These are combined into a global interpolant via a weighted sum, where the weights sum to $1$ everywhere in $\Omega$.

\subsection{Radial basis functions}
\label{ssec:rbf}
First we cover interpolation through collocation with RBFs for some general data set. Then we discuss the somewhat more involved least squares approach, which is the one we will apply. To interpolate over a set of scattered data with a collocation approach, we set each data point to be the centre for an RBF and enforce known function values at these points. This results in a linear system we solve to obtain interpolation weights. Once solved for we can evaluate the interpolant at other desired locations. \\

To be more precise:
Let $s_{f,\Omega}$ be the interpolant of $f\in C^{k}(\Omega)$ on the domain $\Omega\subset \mathbb{R}^{2}$ such that

\begin{equation}
    \label{eq:local_interpol}
    s_{f,\Omega}(\mathbf{y}) = \sum\limits_{j = 1}^{M} \lambda_{j}\varphi_{j}(\mathbf{y}),
\end{equation}
where  $\varphi_{j}(\mathbf{y}) = \varphi(\|\mathbf{y}-\mathbf{x}_{j}\|_{2})$ is the RBF centred at $\mathbf{x}_{j}$ and $\mathbf{y}$ is some point in $\mathbb{R}^{2}$. RBFs are 
 univariate functions from $\mathbb{R}_{+} $ to $ \mathbb{R}$, where  $\mathbb{R}_{+} = \{x\in \mathbb{R} \vert x \geq 0\}$, that disregard geometry. Thus they are mesh free by construction. Given an RBF it forms a basis for the local approximation space $\mathcal{N}_{\varphi}(\Omega)$ referred to in the literature as the \emph{native Hilbert space}. Its nature depends on the RBF of choice, for more details see \cite{Schaback97nativehilbert}. \\

 One option to obtain the unknown \emph{interpolation coefficients} $\lambda_{j}\in\mathbb{R}$  is by collocation at the RBF centres, i.e. for $\mathbf{x}_{i} \in \{\mathbf{x}_{j}\}_{j=1}^{M}$ we require $f(\mathbf{x}_{i}) = \sum\limits_{j = 1}^{M} \lambda_{j}\varphi_{j}(\mathbf{x}_{i})$. This is equivalent to solving 
\begin{equation}
\Phi \Lambda = F
    \label{eq:interpol_linear_sys}
\end{equation}
for $\Lambda = \{\lambda_{j}\}_{j = 1}^{M}$, where $\Phi = \{\varphi_{i}(\mathbf{x}_{j})\}_{i,j = 1}^{M}$  and $F = \{f_{j}\}_{j = 1}^{M}$ with $f_{j} = f(\mathbf{x}_{j})$. The method for obtaining the coefficients through solving \eqref{eq:interpol_linear_sys} with Gaussian--elimination will be referred to as \emph{RBF--Direct}.\\

An RBF $\varphi$ is said to be \emph{positive definite} if the resulting interpolation matrix $\Phi$ is positive definite.
 The matrix $\Phi$ is symmetric due to the nature of RBFs. There is a plethora of positive definite RBFs and we consider only such functions.\\

\begin{figure}[ht]
    \begin{center}
        \includegraphics[width=0.49\textwidth,height=7cm]{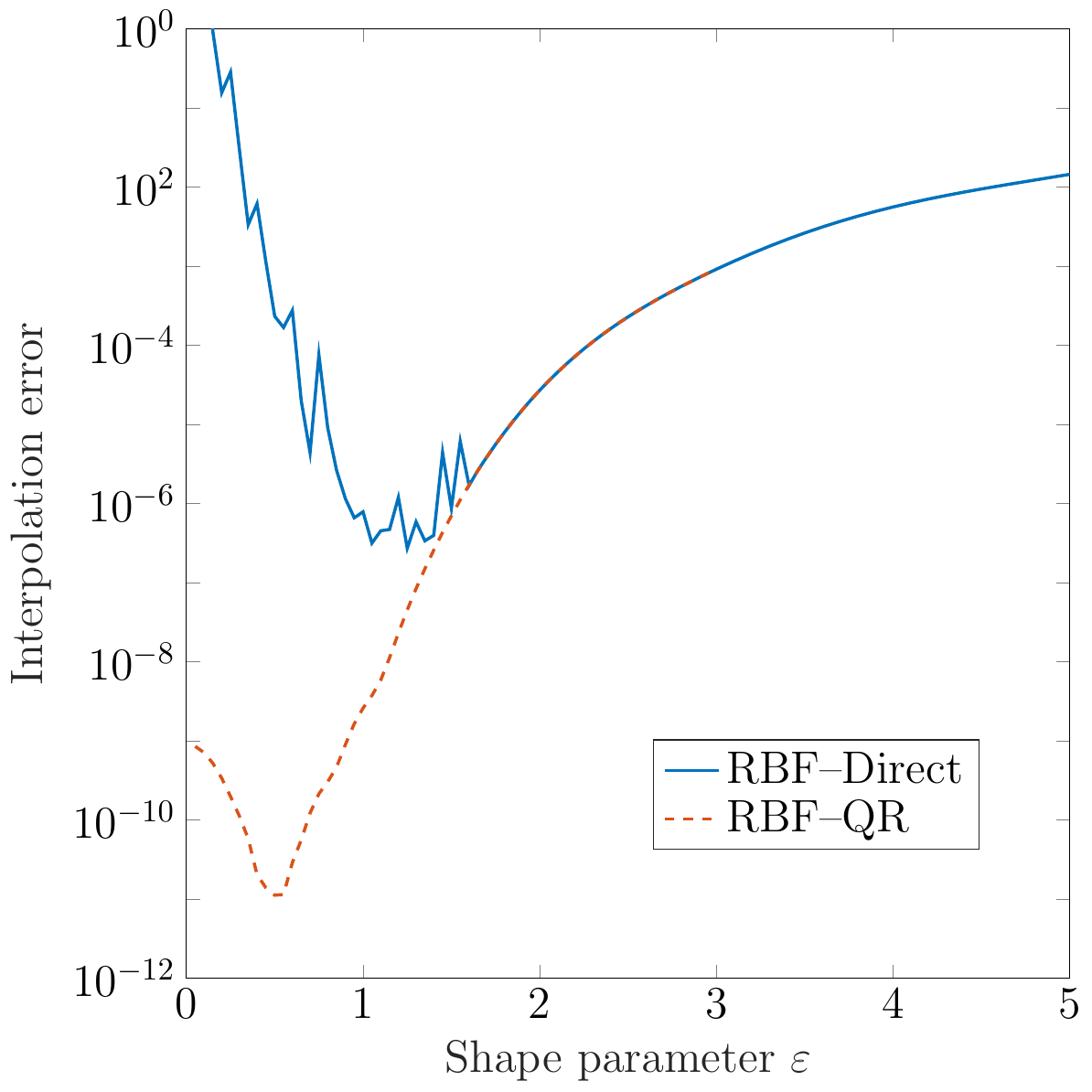}
        \includegraphics[width =0.49\textwidth,height=6.8cm]{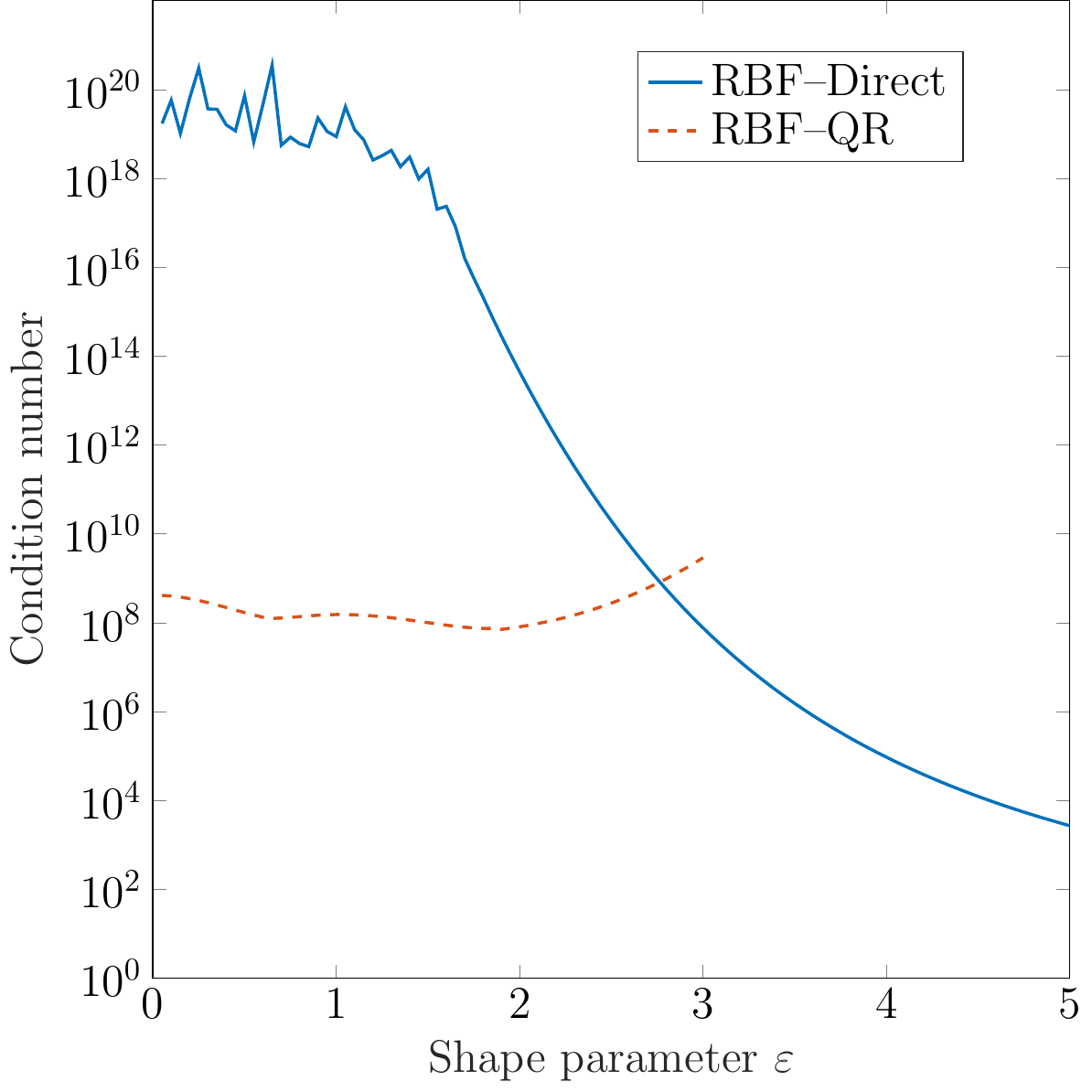}
    \end{center}
    \caption{Illustration of typical interpolation problem for RBF--Direct for $\varepsilon\in (0,5]$ and RBF--QR for $\varepsilon \in (0,3]$. Left: interpolation error. Right: condition number of $\Phi$. }
 \label{fig:increasing_eps}
\end{figure}

 Associated with certain RBFs is the \emph{shape parameter} $\varepsilon\in \mathbb{R}_{+}$, which sets the flatness of the RBF. Let $r = \|\mathbf{y}-\mathbf{x}_{j}\|_{2}$ and consider for example the \emph{Gaussian}
 \begin{equation}
 \varphi(r) = e^{-(\varepsilon r)^{2}},
\label{eq:GaussianRBF}
\end{equation}
for which $\mathcal{N}_{\varphi}(\Omega)$ is a subset in the Sobolev space $W_{2}^{m}(\Omega)$ for any $m$. A function $f$ is an element of $\mathcal{N}_{\varphi}(\Omega)$ for a given $\varepsilon$ if the square root of its Fourier transform decays faster than the Fourier transform of \eqref{eq:GaussianRBF}. For further reading see \citep{Schaback97nativehilbert}. As $\varepsilon$ goes to zero, the RBF becomes increasingly flat and consequently the columns in $\Phi$ become increasingly linearly dependent. Thus the condition number of $\Phi$ grows and the weights $\{\lambda_{j}\}_{j=1}^{M}$ become large and oscillatory. However, these are  purely numerical effects since the RBFs do indeed form the basis for an approximation space with nice properties. In fact, the accuracy of the approximation increases as $\varepsilon $ goes to zero \cite{MR2801193}, but this effect may not be discernible due to high condition numbers, see Figure \ref{fig:increasing_eps}. The appropriate range for $\varepsilon$ is problem dependent and no general optimal value can be chosen, although according to theory the smaller $\varepsilon$ the smaller the error in infinite precision.\\

Apart for Gaussians  other common choices are different compactly supported RBFs, such as the so called \emph{Wu} functions. In the literature  Wu functions are tabulated after their order of regularity at origin, see for example \citep{Fasshauer:2007:MAM:1506263} and Table \ref{tab:Wu}. In our application the regularity at the boundary of the support is of larger interest, where it tends to be greater. This will become apparent in subsequent sections. Thus we will denote the space of continuous functions of regularity $k$
at the origin, but $\tilde{k}$ at the edge of the support, as $C^{k,\tilde{k}}$.\\

\begin{table}
\centering
{\def\arraystretch{1.3}
\begin{tabular}{cc}
\toprule
\textbf{Regularity} & $\varphi(r)$ \\
\midrule
$C^{0,1}$ & $(1- r)_{+}^{2}(2+ r)$ \\
$C^{0,2}$ & $(1- r)_{+}^{3}(8+9 r+3r^{2})$ \\
$C^{2,3}$ & $(1-r)_{+}^{4}(4+16r+12r^{2}+3r^{3})$ \\
$C^{2,4}$ & $(1-r)_{+}^{5}(8+40r+48r^{2}+25r^{3}+5r^{4})$ \\
$C^{4,5}$ & $(1-r)_{+}^{6}(6+36r+82r^{2}+72r^{3}+30r^{4}+5r^{5})$ \\
\bottomrule
\end{tabular}
}
\caption{Wu functions with compact support in $r\in [0,1)$ in $C^{k,\tilde{k}}$: $k$ is the regularity at origin and $\tilde{k}$ is the regularity  at the edge of the support. They are positive definite up to $\mathbb{R}^{3}$. Furthermore $(\cdot)_{+} = \max{(0,\cdot)}$.}
\label{tab:Wu}
\end{table}
 Given a compactly supported RBF of regularity $2k$ at origin the corresponding native Hilbert space is the Sobolev space $W_{2}^{\frac{d}{2}+k+\frac{1}{2}}(\Omega)$, with $\Omega\subset \mathbb{R}^{2}$ \citep{WENDLAND1998258}. If the compact support is small enough they yield a sparse structure for $\Phi$. However this is a trade--off for accuracy, which increases with support size \cite{Fasshauer:2007:MAM:1506263}. We will return to the choice between global and local RBFs.\\

 As mentioned above, the ill--conditioning of $\Phi$ as the shape parameter $\varepsilon$ goes to zero is a purely numerical effect and it can be partially circumvented by replacing RBF--Direct with stable methods, such as \emph{RBF-QR}. This paper does not aim to fully explain the RBF--QR algorithm, instead we refer to \cite{MR2801193} for details. The main idea is that the approximation space spanned by positive definite RBFs contains good solutions to the interpolation problem  \eqref{eq:interpol_linear_sys}, therefore a change of basis can remove the ill-conditioning of $\Phi$. This change of basis is accomplished by the RBF--QR algorithm and its positive effects are clear from Figure \ref{fig:increasing_eps}. Note that for large values of $\varepsilon$ and large distances between RBF centres the RBF--QR is neither needed nor worthwhile, and should hence not be used. \\


As we proceed to the function extension problem, we will have values of $f$ defined on points from a uniform FFT grid. It is a well--known fact in the RBF--community that a uniform point distribution for RBF centres yields a larger condition number. The efficiency of RBF--QR is hampered as well, as it would be limited to collocating at about $100$ uniform grid points. \\

\begin{figure}[!ht]
\centering\includegraphics[width=0.4\linewidth]{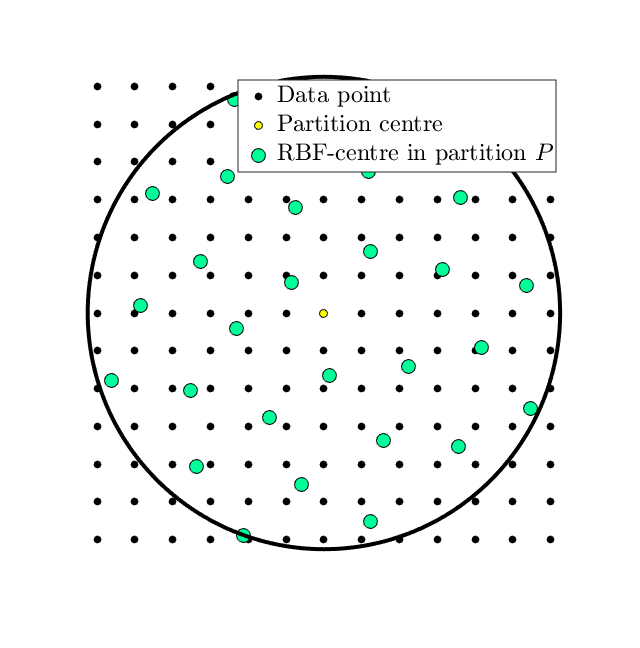}
\caption{Schematic figure of RBF centres and uniformly distributed data points.}
\label{fig:schematic_LeastSquares}
\end{figure}

To obtain a more robust and efficient method, we decouple the centres and the data points, and find an interpolant in the least squares sense, as is done in \citep{Larsson1077950}. Consider a non--uniform distribution of $M$ RBF centres at locations $\myset{\mathbf{x}}{j}{M}$, but the values of $f$ given at $N$ uniform point locations $\myset{\tilde{\mathbf{x}}}{i}{N}$. See Figure \ref{fig:schematic_LeastSquares} for a union of uniform data locations and non--uniform RBF centres. We seek an approximation of $f$ at the non--uniform points, that is $F = \{f(\mathbf{x}_{j})\}_{j=1}^{M}$, in order to create an interpolant with non--uniform RBF centres and data like \eqref{eq:interpol_linear_sys}.  Introduce $\tilde{\Phi} = \{\varphi_{j}(\tilde{\mathbf{x}}_{i})\}_{i,j=1}^{N,M}$ and $\tilde{F} = \{f(\tilde{\mathbf{x}}_{i})\}_{i=1}^{N}$, where  $\tilde{\Phi}$ is an $N\times M$--matrix. Recall that $\Phi$ and $\Lambda$ is of size $M\times M$ and $M\times 1$, respectively. Then, replacing $\Lambda$
by $\Phi^{-1}F$  from \eqref{eq:interpol_linear_sys}, we have
\begin{equation}
\tilde{\Phi} \Lambda = \tilde{F} \Leftrightarrow \tilde{\Phi} \Phi^{-1} F = \tilde{F},
\label{eq:LS}
\end{equation} 
and let $A = \tilde{\Phi} \Phi^{-1}$. Obtaining the unknown $F$ from \eqref{eq:LS} is a least squares problem and numerically we require that $N$ is sufficiently larger than $M$ to yield good solutions. This approach allows us to use non--uniform RBF centres which significantly improves the stability, but still lets the data be represented on the uniform grid. Furthermore, the ill--conditioning of $\Phi$, associated with the shape parameter $\varepsilon$, is reduced by the use of RBF--QR. The RBF--QR algorithm is intended for a formulation as \eqref{eq:LS}, since it computes $A$, rather than $\Phi^{-1}$, which acts as a mapping of data from non--uniformly  to uniformly distributed locations. With $F$ known we can interpolate to obtain values $F_{\mathbf{y}}$ at some locations $\myset{\mathbf{y}}{k}{K}$ by creating  $A_{\mathbf{y}} = \Phi_{\mathbf{y}}\Phi^{-1} $ where $\Phi_{\mathbf{y}} = \{\varphi_{j}(\mathbf{y}_{k})\}_{k,j=1}^{K,M}$; then $F_{\mathbf{y}} = A_{\mathbf{y}}F$ which we set to be $s_{f,\Omega}(\mathbf{y}_{k})$ for $k = 1,\ldots,K$.\\

The distribution of non--uniform RBF centres is the same as in \cite{Larsson1077950} and is a \emph{Vogel node} distribution defined as
\begin{equation}
\mathbf{x_{j}} = \sqrt{\frac{j}{M}}\left(\cos{(j\,\pi(3-\sqrt{5}))},\,\sin{(j\pi(3-\sqrt{5}))}\right),\quad j = 1,\ldots,M,
\label{eq:Vogelnodes}
\end{equation}
 and is quasi--uniform. In Figure \ref{fig:schematic_LeastSquares} the $28$ mint green dots are Vogel nodes. Such a distribution of RBF centres is near--optimal and the RBF--QR algorithm performs well up to about $400$ centres. Recall that in order obtain good approximations by solving the least squares system $N$ must be sufficiently larger than $M$.\\

A drawback when using Gaussians is their global nature and consequently the resulting interpolation matrix $\Phi$ is large and dense. An alternative would be to use compactly supported RBFs, but to obtain good approximations their support needs to be large, resulting again in a non--sparse matrix structure. This encourages the implementation of a partition of unity method, which is often used in combination with RBFs. It decouples the size of $\Phi$ and the choice of the RBF used for interpolation.

\subsection{Partition of Unity}
\label{ss:pum}
The idea of a partition of unity approach is to combine local approximations $s_{f,\Omega_{i}}$ of the function $f$ on partitions $\Omega_{i}$, $i = 1,\ldots,N_{p}$. These partitions form the set $\{\Omega_{i}\}_{i = 1}^{N_{p}}$. They are overlapping and constitute a covering of $\Omega$, meaning 
\begin{equation}
\label{eq:covering}
\bar{\Omega}\subset \bigcup\limits_{i=1}^{N_{p}}\bar{\Omega}_{i}. 
\end{equation}
Associated with this covering we construct a family of compactly supported and continuous functions $\{w_{i}\}_{i = 1}^{N_{p}}$ such that $\text{supp}(w_{i})= \Omega_{i}$ for every $i = 1,\ldots,N_{p}$ and 

\begin{equation}
\label{eq:pou1}
\sum\limits_{i = 1}^{N_{p}} w_{i}(\mathbf{y}) \equiv 1, \quad \forall \mathbf{y}\in \bigcup\limits_{i=1}^{N_{p}}\bar{\Omega}_{i}.
\end{equation}
Each weight function $w_{i}$ corresponds to a partition $\Omega_{i}$ in the covering. They are constructed as a weighted average of compactly supported RBFs which we denote $\psi$, to distinguish them from the Gaussians $\varphi$. This is often referred to in the literature as Shepard's method \cite{shepard}:

\begin{equation}
    w_{i}(\mathbf{y}) = \frac{\psi_{i}(\mathbf{y})}{\sum\limits_{j = 1}^{N_{p}}\psi_{j} (\mathbf{y})},
    \label{eq:weightfunction}
\end{equation}
where $i\in{1,2,\ldots,N_{p}}$. An RBF  $\psi_{i}$ is centred at $\mathbf{p}_{i} \in \mathbb{R}^{2}$ and its support defines the associated partition $\Omega_{i}$. The set $\{\mathbf{p}_{i}\}_{i=1}^{N_{p}}$ contains all partition centres.\\

 We will use the same distribution of non--uniform RBF centres for $\{\varphi_{j}\}_{j=1}^{M}$ for every partition. Thus the centres $\{\mathbf{x}_{j}\}_{j=1}^{M}$ can then be selected for near--optimal approximation properties, where Vogel node distribution is one such example, see \eqref{eq:Vogelnodes}. The approximation $s_{f}$ of $f$ on $\Omega$ is a weighted sum of the local approximations $s_{f,\Omega_{i}}$, i.e.
\begin{equation}
\label{eq:rbfpum}
    s_{f}(\mathbf{y}) = \sum\limits_{i = 1}^{N_{p}}w_{i}(\mathbf{y})s_{f,\Omega_i}(\mathbf{y})
\end{equation}
where the local approximations are obtained through interpolation with RBFs:

\begin{equation}
    s_{f,\Omega_i}(\mathbf{y}) = \sum\limits_{j=1}^{M} \lambda_{j}\varphi_{j}(\mathbf{y}),
\end{equation}
 see subsection \ref{ssec:rbf}. From now on we will always use Gaussians as RBFs for interpolation and  all partitions are discs. As we proceed with function extension the local approximations are only required for partitions partially in $\Omega$. For partitions entirely in $\Omega$ constructing  a local extension is redundant, since $f$ is known there.


\section{Function extension by PUX} 
\label{s:functionextesion}

In this section we describe how to construct a compactly supported  function extension of high global regularity. To make the description of the PUX--method accessible we will first  extend a function with extrapolation, using tools presented in the previous section. For this extension we do not control how $f^{e}$ approaches zero. Thereafter we apply the full PUX--method to construct a compactly supported function extension with a chosen regularity.

\subsection{Extrapolation}
\label{ss:Extrapolation}
\begin{figure}[ht]
    \begin{center}

\begin{tikzpicture}
\node (bad) {\includegraphics[trim=2.6cm 1.7cm 1.7cm 1.3cm, clip, width=0.3\textwidth]{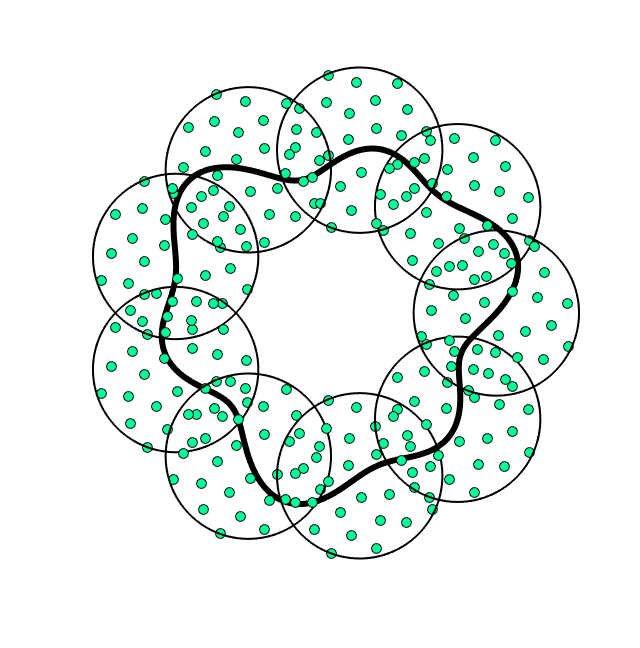}};
\node (good) [right=of bad] {\includegraphics[trim=2.2cm 2.4cm 2cm 1.6cm, clip,  width =0.3\textwidth]{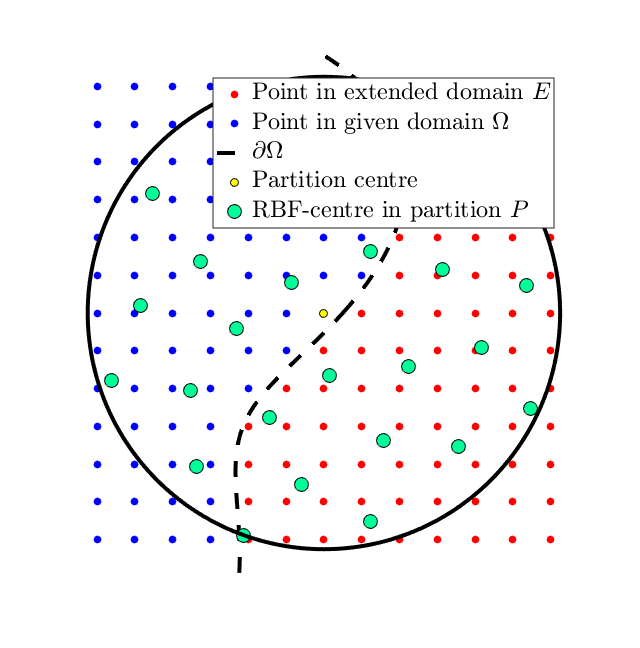}};

\draw [ultra thick,black,->] (2.2,0.7) to[bend left] (3.6,1);
\end{tikzpicture}

    \end{center}
\caption{Left: Schematic figure of distribution of extension partitions along $\partial \Omega$ for a complex domain. The green markers correspond to RBF--centres and the distribution is repeated for every partition. Right: Schematic figure of identifying uniformly distributed points as inside or outside $\Omega$.}
\label{fig:ls_entire}

\end{figure}

For extrapolation of $f$ defined on a bounded and closed domain $\bar{\Omega}\subset \mathbb{R}^{2}$ the following procedure is used. Partition centres $\{\mathbf{p}_{i}\}_{i=1}^{N_{p}}$, i.e. centres for the compactly supported RBFs we denote $\psi_{i}$, are distributed uniformly with respect to arc length along the boundary $\partial \Omega$ as in the left image in Figure \ref{fig:ls_entire}.  Each partition centre is associated with a partition for which a local extension is constructed. They are referred to as \emph{extension partitions}. The partition centres are moved to the closest uniform grid point in $\Omega$, such that if all extension  partitions have the same radius $R_{p}$, a single matrix $A$ from expression \eqref{eq:LS} can be precomputed and reused for every partition $\Omega_{i}$ to approximate the local interpolants $s_{f,\Omega_{i}}$, $i=1,\ldots, N_{p}$. Doing so is vital for the efficiency for PUX, as RBF--QR is required for accuracy but is $10$ to $15$ times more computationally expensive than RBF--Direct.\\

We introduce some useful notation: let $\mathbf{y}_{E} = \{\mathbf{y}\in E\}$ and $\mathbf{y}_{\Omega} = \{\mathbf{y}\in \Omega\}$, where $E$ is the complement of $\bar{\Omega}$ relative $B$. Any combination of subindices simply means points in the associated intersection, e.g. $\mathbf{y}_{E,i} = \{\mathbf{y}_{E}\in\Omega_{i}\}$. The point distribution between partitions differ only with regard to which points that belong to $\Omega$ or  $E$. See Figure \ref{fig:ls_entire} and compare with Figure \ref{fig:schematic_LeastSquares}; observe that we have RBF centres outside of $\Omega$. For some partition centre $\mathbf{p}_{i}$ we precompute a single matrix $A$ based on all uniform points and Vogel nodes within distance $R_p$ of the centre. For a given partition $\Omega_{i}$ we identify the uniform points $\mathbf{y}_{\Omega_{i}}$ as belonging to $\Omega$ or to $E$. By simple row manipulation we rewrite \eqref{eq:LS} as
 \begin{equation}
 \label{eq:interpol_LS_sys}
    \begin{bmatrix}
    A_{\Omega}\\
   A_{E}
\end{bmatrix}
F=
\begin{bmatrix}
\tilde{F}_{\Omega} \\
\tilde{F}_{E}
\end{bmatrix}
 \end{equation}
where the subscripts denote which set the points $\mathbf{y}_{\Omega_{i}}$ belong to. Since $\tilde{F}_{\Omega}$ is known we can solve $A_{\Omega}F = \tilde{F}_{\Omega}$ for $F$ in least squares sense. Thus the number of defined interior values for $f$, i.e. the size of $\tilde{F}_{\Omega}$, must be greater than the number of RBF--centres, i.e. Vogel nodes,  in order to have an overdetermined least squares problem. The \emph{local extension} is computed as
\begin{equation}
s_{f,\Omega_{i}}(\mathbf{y}_{\Omega_{i},E}) = \tilde{F}_{E} = A_{E}F.
\label{eq:extrapolatingf}
\end{equation}
To ease notation we refer to $s_{f,\Omega_{i}}(\mathbf{y}_{\Omega_{i},E})$ as $f^{e}_{i}$ when no ambiguity can arise. The global extension $f^{e}$ is the weighted combination of all local extension evaluated for $\cup_{i}\mathbf{y}_{\Omega_{i}, E}$, but for $\cup_{i}\mathbf{y}_{\Omega_{i}, \Omega}$ we set $f^{e} = f$ and $f^{e}=0$ otherwise. The extension has compact support, but no means are taken to control how it decays. Hence in most events the extension will be discontinuous over the boundary of its support. A remedy is presented in the following subsection.

\begin{figure}[ht] 
  \begin{subfigure}[b]{0.5\linewidth}
    \centering
    \includegraphics[trim=2cm 0cm 2cm 0cm, clip,width=\linewidth]{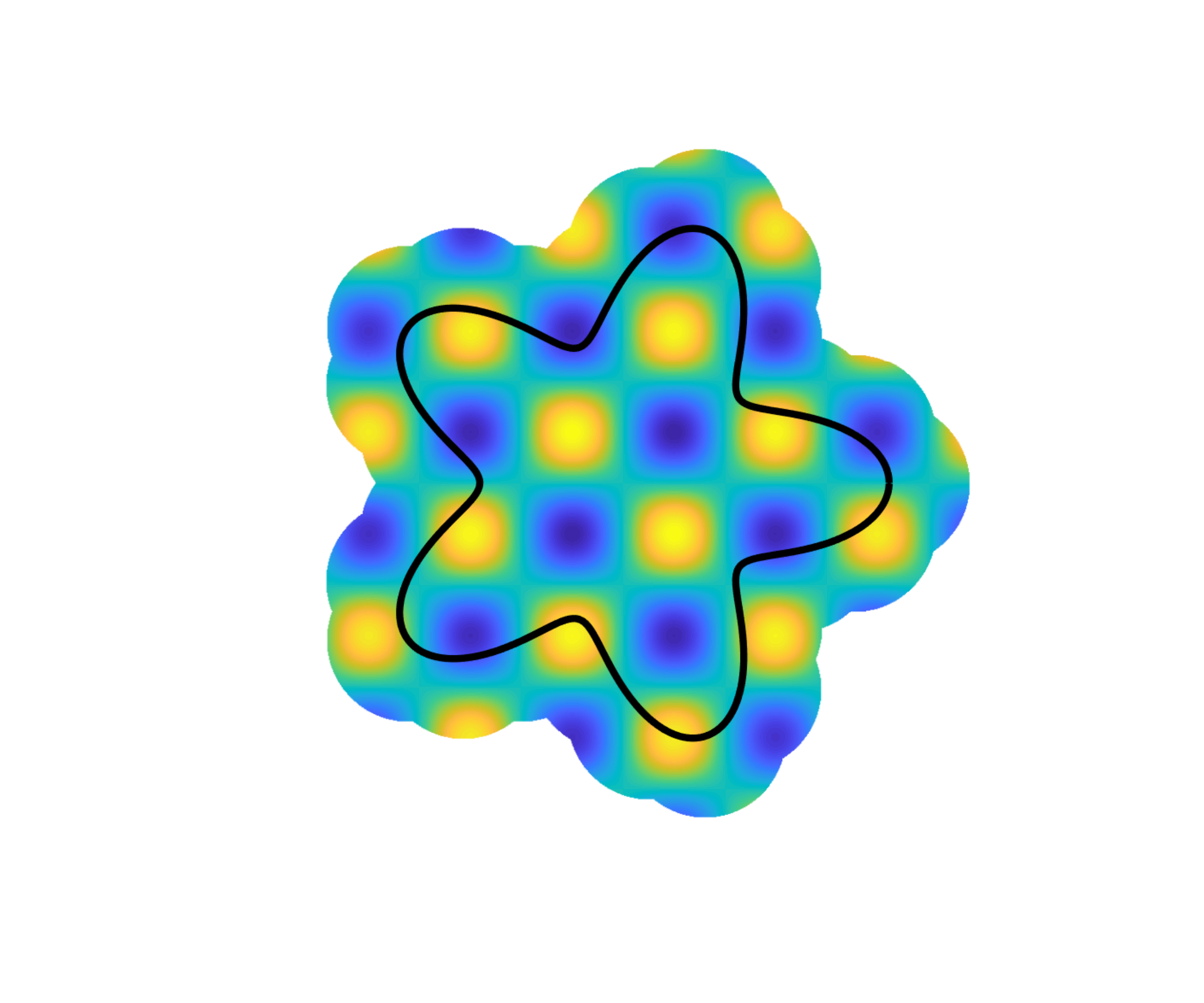} 
    \caption{Function $f$ from \eqref{eq:results1} on starfish shaped domain $\Omega$.} 
    \label{fig:puxA} 
    \vspace{4ex}
  \end{subfigure}
  \begin{subfigure}[b]{0.5\linewidth}
    \centering
    \includegraphics[trim=2cm 0cm 2cm 0cm, clip,width=\linewidth]{images/f_extrapol-eps-converted-to.pdf} 
    \caption{Extrapolation of \eqref{eq:results1} as in subsection \ref{ss:Extrapolation}.} 
    \label{fig:puxB} 
    \vspace{4ex}
  \end{subfigure} 
  \begin{subfigure}[b]{0.5\linewidth}
    \centering
    \includegraphics[trim=2cm 0cm 2cm 0cm, clip,width=\linewidth]{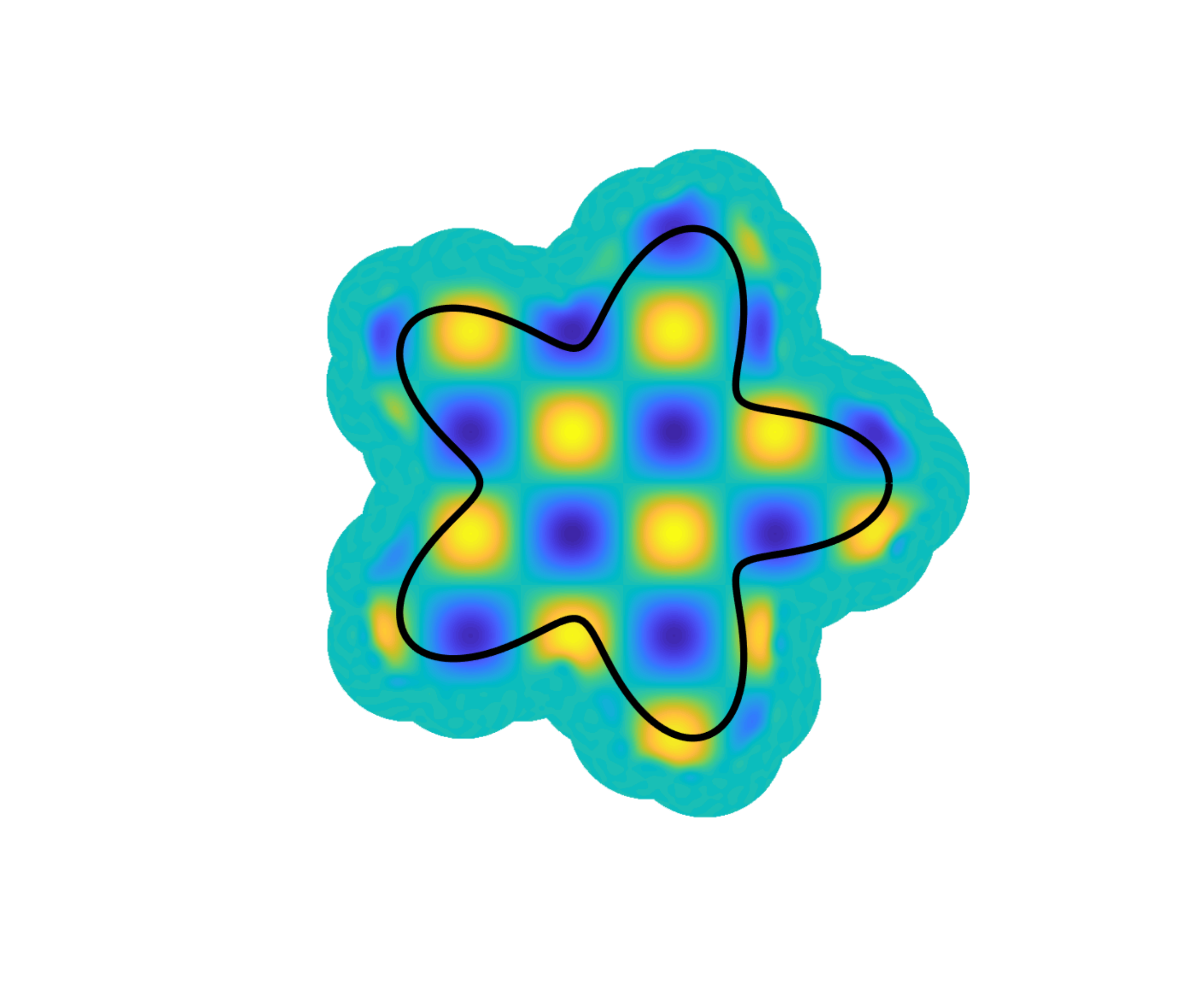} 
    \caption{Blended with a layer of\\ zero--partitions, see subsection \ref{ss:CSextension}.} 
    \label{fig:puxC} 
  \end{subfigure}
  \begin{subfigure}[b]{0.5\linewidth}
    \centering
    \includegraphics[trim=2cm 0cm 2cm 0cm, clip,width=\linewidth]{images/f_supress-eps-converted-to.pdf} 
    \caption{High regularity function extension $f^{e}$ with compact support, see \eqref{eq:functionextension}.} 
    \label{fig:puxD} 
  \end{subfigure} 
  \caption{Illustration of various steps of function extension with PUX.}
  \label{fig7} 
\end{figure}

\FloatBarrier









\begin{figure}[ht]
    \begin{center}

\begin{tikzpicture}
\node (bad) {\includegraphics[trim=2.4cm 1.7cm 1.4cm 1cm, clip, width=0.3\textwidth]{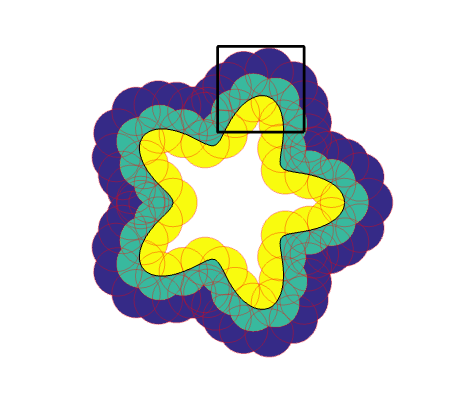}};
\node (good) [right=of bad] {\includegraphics[trim=1.8cm 1.2cm 1.2cm 0.8cm, clip, width =0.3\textwidth]{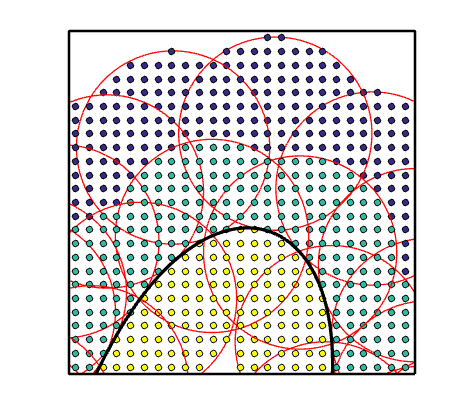}};
\draw [ultra thick,black,->] (1.2,2.5) to[bend left] (3.6,2.5);
\end{tikzpicture}

    \end{center}
\caption{Schematic image for function extension from a star shaped domain $\Omega$ given by the black border. Observe that in this figure the partitions are not centred at uniform grid points. The red overlapping circles are the partitions. The yellow section corresponds to the uniform data  points used for creating the local extension $f^{e}_{i}$, the blue section to points where $f^{e} = 0$ and the green sector is a blend of the two.}
\label{fig:schematic_extension_supression}

\end{figure}

\subsection{Compactly supported function extension}
\label{ss:CSextension}
To construct a compactly supported function extension $f^{e}$ of $f$, which is continuous or of higher regularity as it is extended by zero outside of its support, we modify the extension obtained in the previous subsection. Yet another layer of partitions is added such that it overlaps the extension partitions. These partitions do not intersect $\bar{\Omega}$ and no interpolation is performed  for them. They are referred to as \emph{zero partitions} and the corresponding set is denoted as $\myset{\Omega^{0}}{i}{N^{0}_{p}}$. The associated $s_{f,\,\Omega^{0}_{i}}$ is set to be identically equal to zero for $i = 1,\ldots, N^{0}_{p}$.\\

The weight function in a partition, as defined in \eqref{eq:weightfunction}, is zero at the boundary of the partition. Hence, as the local extension in the first layer of partitions are blended with the zero values in the zero partitions, the global extension that is defined will be forced to zero over the overlapping region. See Figure \ref{fig:puxC} where the local extensions are suppressed to zero in the overlaps. Hence zero partitions should be placed such that $f^{e}$ has a controlled decay to zero and that the size of the overlap with extension partitions are about the same. Thus the global extension will in these parts have the same regularity as $w$, i.e. the compactly supported RBF $\psi$. Observe that there always is jump over $\partial \Omega$, since an analytical expression is used for points inside $\Omega$ and an approximation outside. It is not numerically  discernible if the approximation is good enough. However, if poor then $f^{e}$ is an extension of a poor approximation of $f$, i.e. another function. This occurs if e.g. the number of RBF--centres per partition is too few. Consequently the jump $\partial \Omega$ will be of such magnitude that $f^{e}$ will behave as a discontinuous function.\\

Once $f^e$ is obtained it is extended by zero to a box $B = [-L,L]^{2}$ that embeds the support of $f^e$. The expression for a function extension obtained by PUX is 
\begin{equation}
  \label{eq:functionextension}
      f^{e}(\mathbf{y}) = \begin{cases}f(\mathbf{y}),\quad \mathbf{y}\in\Omega,\\\sum\limits_{i = 1}^{N_{p}+N_{p}^{0}}w_{i}(\mathbf{y})s_{f,\Omega_i}(\mathbf{y}),\quad \mathbf{y} \in \bigcup\limits_{i=1}^{N_{p}}\mathbf{y}_{\Omega_{i}, E},\\
      0, \quad \text{otherwise}.
      \end{cases}
\end{equation}

The sum also includes the zero partitions to emphasise that the weight functions in extension partitions are affected by them; even though the local values in a zero partition are identically zero. The function $f^{e}$ evaluated on a uniform grid in $B$ can be  used to solve
 \eqref{eq:poissonR2}--\eqref{eq:poissonR2L2} for the particular solution to the full problem. \\

\section{Sources of numerical errors}
\label{s:errorAnalysis}

Given a  function $f$ defined on a bounded domain $\Omega$, there is no unique compactly supported extension to $\mathbb{R}^{2}$. To measure the quality of the function extension obtained by PUX for a set of parameters, we study the numerical solution to the full problem \eqref{eq:poisson}--\eqref{eq:poissonBC}. Consequently, the errors directly associated with PUX  cannot be isolated and analysed separately, but we can give an account for sources of numerical errors. In this section we discuss the local interpolation error, the error originating from approximating local interpolants with least squares, errors dependent on the choice of the compactly supported RBF to construct the weight function, errors from using an FFT--solver for the free--space Poisson equation \eqref{eq:poissonR2} with the truncated Green's function and errors associated with boundary integral method for solving the Laplace equation \eqref{eq:laplace_mod}--\eqref{eq:laplace_modbc}.\\

First we briefly discuss the error from solving the Laplace equation. Assume $\partial \Omega$ is sufficiently resolved with $16$ point Gauss--Legendre panels. Then the error associated with the boundary integral method for solving the Laplace equation converges rapidly as the number of Gauss-Legendre panels increases. The resolution of $\partial \Omega$ is chosen a priori based on $g$, therefore it is important that the modified boundary conditions $g-u^{p}$ are not much harder to resolve than $g$. If this is true, then the error associated with the boundary integral method can with efficiency and without difficulty be controlled such that it does not dominate.\\


Now we turn our attention to the errors connected to function extension and solving the free--space Poisson equation \eqref{eq:poissonR2}. We  start with the local interpolation error in some partition $\Omega_{i}$; to initially simplify the analysis assume  $\{\lambda_{j}\}$ in \eqref{eq:local_interpol} is obtained by collocation, i.e. not with least squares. If the RBFs $\{\varphi_{j}\}$ are Gaussians, then according to \cite{Rieger2009} the following estimate for the interpolation error holds:

\begin{equation}
\|f - s_{f,\Omega_{i}}\|_{L_{\infty}(\bar{\Omega}_{i})} \leq 3 e^{C \log{(h)}/\sqrt{h}}\,\|f\|_{\mathcal{N}_{\varphi}},
\label{eq:errorEst}
\end{equation}
where $C$ is a constant depending on the space dimension, the norm and $\varphi$, but not on $f$ or $h$. Here $h$ is the \emph{fill distance}: the diameter of the largest ball possible to fit between the RBF centres in a partition and thus a measure of the density of the RBF--centres. Further $\|\cdot\|_{\mathcal{N}_{\varphi}}$ is the norm associated with the native space $\mathcal{N}_{\varphi}$. Although technically only true for $f\in \mathcal{N}_{\varphi}$\, it works well in practise for smooth functions and especially bandlimited functions. For further discussion see \cite{Larsson2005}. From \eqref{eq:errorEst} we conclude that the interpolation error in each partition for collocation converges spectrally as we increase the number of RBF--centres within a partition. We note that the estimate \eqref{eq:errorEst} still holds when $s_{f,\,\Omega_{i}}$ is obtained with least squares, since the same RBFs are used and therefore the approximation space is the same as for collocation \cite{Larsson1077950}. However, this will not be true for the extension.\\

To understand what happens with the interpolation error when the local interpolants are weighted together with a partition of unity method we need the concept of \emph{regular covering}: Assume that each $\mathbf{y}\in \Omega$ only belongs to a finite number of $\Omega_{i}$, that each partition contains enough data points to allow a unique interpolant and that each partition $\Omega_{i}$ satisfies an interior cone--condition \cite{Wendland02fastevaluation}. A covering $\{\Omega_{i}\}_{i = 1}^{N_{p}}$ satisfying these requirements is denoted a regular covering. In our applications, these conditions are easily met and we will assume that our coverings are regular.\\

A remarkable property of partition of unity is that the global approximation order of $s_{f}$ inherits the local one for $s_{f, \Omega_{i}}$, under the assumptions that the covering $\{\Omega_{i}\}_{i = 1}^{N_{p}}$ is regular \cite{Wendland02fastevaluation}. This property can be understood intuitively, just consider the following example. Let $\mathbf{y}\in \Omega_{1}\cap \Omega_{2}$ and assume the local approximant is exact up to some tolerance $\epsilon_{i}$, i.e. $\|f-s_{f,\Omega_{i}}\|_{L_{\infty}(\bar{\Omega}_{i})}\leq \epsilon_{i}$ for $i = 1,2$. Then
\begin{equation}
\label{eq:rbfpumerror}
    |f(\mathbf{y})-s_{f}(\mathbf{y})| \leq \sum\limits_{i = 1}^{2}  w_{i}(\mathbf{y})|f(\mathbf{y})-s_{f,\Omega_{i}}(\mathbf{y})| \leq \underbrace{(w_{1}(\mathbf{y}) + w_{2}(\mathbf{y}))}_{=1}\max(\epsilon_{1},\epsilon_{2})
\end{equation}
and the potential error from low order weight functions is not noticeable, as they form a partition of unity and thus sum to one at every point belonging to a partition. Instead,  the global interpolant inherits the local approximation order.\\


Now we consider the function extension and discuss its regularity and influence on the error, still assuming that the local interpolants are obtained through collocation instead of approximated with least squares. The speed at which the extension goes to zero is controlled by the choice of RBF we use as $\psi$ and on the size of the overlapping region. A smaller region means a more rapid decay, thus a high grid resolution might be required to resolve the extension. Therefore it is preferable, from this point of view, to have as much overlap between extension partitions and zero--partitions as possible, and to have a large partition radius. On the other hand, the extension may not be well--behaved far from the boundary and a greater partition radius requires more RBF-centres to resolve $f$. The RBF--QR algorithm is not capable of removing ill--conditioning, associated with small values of $\varepsilon$, for $A$ with over $400$ RBF centres. Furthermore, the larger the partition radius $R_{p}$ is the more we lose locality.\\
 
 Increasing the number of RBF centres does not necessarily increase the accuracy, instead it can make the situation worse. Compare with polynomial interpolation, where a higher order polynomial basis means more oscillations. A similar effect is present for RBF--interpolation. A remedy is to make the partitions smaller, instead of increasing the number of RBF centres. However, this requires a finer resolution of the uniform grid, as smaller partitions implies a shorter span for the extension to go to zero.\\

Theoretically, things get more complicated when we consider the function extension obtained by least squares. With collocation the interpolant and $f$ would agree on all uniform data locations in $\Omega$, where as for the least squares interpolant, there can be a discrepancy at these points. Defining $f^{e}$ according to \eqref{eq:functionextension}, where the original values of $f$ are used inside $\Omega$, there will hence be a discontinuity across $\partial \Omega$ of the size of the error in the least squares interpolant. However, this discontinuity is a technicality, since it can be made arbitrary small by controlling \eqref{eq:errorEst}. If no such measures are taken the error for solving the full problem decays with second order as the uniform grid is refined. Since $A$ can be precomputed once and used for all partitions, we can set the number of RBF--centres $M$ sufficiently large and make the partitions appropriately small in order not to suffer from a poor approximation of the interpolant due to least squares. 
\\

We now discuss the weight functions' influence on the accuracy for numerically solving \eqref{eq:poisson}--\eqref{eq:poissonBC} . The intersection between a zero and an extension partition is the distance over which the weight functions suppress the extension to zero. This implies that a high grid resolution may be required not to resolve $f$ on $\Omega$, but to resolve the suppression to zero over the aforementioned intersection. Further the extension will inherit the regularity of the compactly supported RBF used to create the weights $w$. Thus we may expect an asymptotic convergence of $4+\tilde{k}$ for an RBF with regularity $C^{k,\,\tilde{k}}$: $2$ orders for solving the free--space Poisson equation plus $2+\tilde{k}$ since the $\tilde{k}$:th derivative has bounded variation \cite{spectralMethodsInMatlab}. The greater the regularity of the Wu RBF the harder it is to resolve, thus requiring a finer uniform grid. Therefore we will not coinsider infinitely smooth compactly supported RBFs as weight functions, such as the construction by  Ying et al. \citep{Ying2006} . In section \ref{s:numericalResults} we show how to optimally choose the regularity $\tilde{k}$ of the compactly supported RBF, given a resolution of the uniform grid on $B$ and partition radius $R_{p}$.
\\

Finally we briefly mention errors associated with evaluating the solution to the free--space Poisson equation with a truncated kernel by FFTs and for obtaining $u^{P}$ on $\partial \Omega$ with non--uniform FFT. Replacing the Green's function $K$ with the truncated Green's function $\tilde{K}$ in \eqref{eq:poissonPartSol} involves no approximation. The main source of error is resolving the Fourier transform of $\tilde{K}$ and $f^{e}$ with FFTs. Given a smooth right hand side $f^{e}$ the error for solving \eqref{eq:poissonR2} decays spectrally as the uniform grid is refined, assuming method parameters are set appropriately \cite{MR3535728}. Concerning the non--uniform FFT it takes a given tolerance as input argument, and assuming sufficiently many Fourier coefficients are available no greater error than the set tolerance  will be introduced.\\

 \section{Summary and implementation of entire solution procedure}
\label{s:summary}
Guidelines for picking appropriate values for the numerical implementation are given in the next section. In this section the computational procedure for solving the full problem \eqref{eq:poisson}--\eqref{eq:poissonBC} is summarised, alongside implementation details:
\begin{itemize}

\item Discretise $\partial \Omega$ into $N_{\partial \Omega}$ Gauss--Legendre panels, each with $16$ Gauss--Legendre points. Set $R_{p}$ as radius for all extension partitions and let them overlap by slightly more than a radius. This yields a number $N_{p}$ of extension partitions  $\{\Omega_{i}\}_{i=1}^{N_{p}}$, with corresponding partitions centres $\{\mathbf{p}_{i}\}_{i=1}^{N_{p}}$, distributed uniformly with respect to arc length along $\partial \Omega$. \\

\item Let $\Omega \subset B = [-L,L]^{2}$, which must contain all $\mathbf{y}\in\cup_{i}\bar{\Omega}_{i}$, and let $E$ be the complement of $\bar{\Omega}$ relative $B$. At the locations  $E\backslash \cup_{i}\bar{\Omega}_{i}$ the extended function will be zero and need not to be included, see \eqref{eq:functionextension} or Figure \ref{fig:results1covering}. For some $N_{u}$  construct a uniform grid on $B$ with resolution $2L/N_{u}$ to be used by the FFT--solver. Sort the $N_{u}^{2}$ uniform grid points as \emph{inside} or \emph{outside} $\bar{\Omega}$. To identify a point $\mathbf{y}$ as belonging to $\Omega$ or $E$, one can evaluate \eqref{eq:doublelayer_mu1} for it.
 \\

\item Relocate each partition centre $\mathbf{p}_{i}$ to the closest uniform grid point inside $\Omega$ to allow for precomputation of $A$ and to avoid reducing the regularity of $f^{e}$ by evaluating the RBFs at their origin.

\item Pick a compactly supported RBF for $\psi$ to construct the weight functions $\{w_{i}\}_{i=1}^{N_{p}+N^{0}_{p}}$ with Shepard's method \eqref{eq:weightfunction}. We use Wu--functions, see Table \ref{tab:Wu} and choose among them based on $\tilde{k}$, which denotes the compactly supported RBF's regularity, with a neighbourhood around the origin excluded. This value sets the regularity of the global extension $f^{e}$.

\item The shape parameter $\varepsilon$ for the Gaussians used as the RBF for interpolation \eqref{eq:local_interpol} needs to be set. The number of RBF centres $M$, i.e. Gaussians, inside each partition and which distribution must be specified too. We use a Vogel node distribution \eqref{eq:Vogelnodes} and for a given mesh size $N_{u}$ and partition radius $R_{p}$ the amount $M$ must be set such that \eqref{eq:LS} is an overdetermined system. If the partition radius $R_p$ is appropriately small in relation to the curvature of $\partial \Omega$, each partition will contain roughly the same amount of uniform data locations where $f$ is known. This does not only give the same interpolation qualities on each patch, but is a measure to ensure that for each patch the least squares problem \eqref{eq:LS} is sufficiently overdetermined. In the next section we will give a guide to choosing all aforementioned parameters. 

\item To precompute the matrix $A$ pick some point on the uniform grid.  Within distance $R_{p}$ of this grid point, find the uniform grid points  $\{\tilde{\mathbf{x}}_{i}\}_{i=1}^{N}$  and distribute a set of Vogel nodes $\{\mathbf{x}_{j}\}_{j=1}^{M}$, see Figure \ref{fig:schematic_LeastSquares}.  For these points we compute $A = \tilde{\Phi} \Phi^{-1}$ with RBF--Direct, where $\tilde{\Phi} = \{\varphi_{j}(\tilde{\mathbf{x}}_{i})\}_{i,j=1}^{N,M}$ and $\Phi = \{\varphi_{j}(\mathbf{x}_{i})\}_{i,j=1}^{M,M}$, see paragraph above expression \eqref{eq:LS}. If the condition number of $\Phi$ is not of moderate size, then $A$ is recomputed with RBF--QR. An open--source implementation in \textsc{Matlab} for computing $A$ with RBF--QR can be found at \cite{RBFQRlink}. This is done only once, and $A$ is reused for all extension partitions.\\

\item  For each partition $\Omega_{i}$ the precomputed  matrix $A$ needs to be separated into $A_{\Omega}$ and $A_{E}$, as in \eqref{eq:interpol_LS_sys}. The patch specific separation depends on which uniform points $\mathbf{y}_{\Omega_{i}}$ that belong to $E$ or $\Omega$. Standard solvers are used to solve the least squares systems, see the paragraph below expression \eqref{eq:interpol_LS_sys}, in order to evaluate the local extension $f^{e}_{i}$ by \eqref{eq:extrapolatingf}.\\


\item It is also possible to precompute the local component $\psi$ used for constructing the weight functions \eqref{eq:weightfunction}. As for precomputing $A$, simply evaluate $\psi$ for all uniform points within distance $R_p$ of some partition centre. These values can be used for all partitions to construct the partition of unity weights.\\

\item The $N_{p}^{0}$ zero partitions must at least overlap the boundary of the union of $\Omega$ and the extension partitions in order to have a controlled decay of $f^{e}$ to zero. One simple way of achieving this is to distribute $N_{p}^{0}= 2 N_{p}$ points uniformly, with respect to arc length, on the boundary. Place the zero--partitions in the normal direction at a distance of $R_{p}$ from the boundary $\partial \Omega$. Initially their radius is set to $R_{p}$, but as they should not intersect $\bar{\Omega}$ rescaling may need to be required. Thus the zero--partitions can have varying radius.\\

\item Combine the  local extensions $\{f^{e}_{i}\}$ by weight functions to obtain the global extension $f^{e}$ of $f$, given by \eqref{eq:functionextension}, for each point from the uniform grid in  $B$.\\

\item Since $f^{e}$ has compact support in a box $B$ embedding $\Omega$ we can apply the method in \cite{MR3535728}, as explained in subsection \ref{ss:homogeneuous}. Hence we use the precomputed truncated spectral representations of the Green's function \eqref{eq:truncatedSpectralRepofGreen}, to evaluate the solution to \eqref{eq:poissonR2}--\eqref{eq:poissonR2L2} with an FFT, upsampled by a factor of $2$. This yields the particular solution $u^{p}$ on a uniform grid in $B$ and also the Fourier coefficients of $u^{P}$. Evaluate $u^{P}$ at the Gauss--Legendre points given by the Gauss--Legendre panels along $\partial \Omega$ with a non--uniform FFT. The non--uniform FFT we apply is described in \cite{MR2115056} and we use their open source library. As input it takes an error tolerance which we set to $10^{-14}$. These values are used to obtain the modified boundary conditions $g - u^{p}$ for \eqref{eq:laplace_mod}--\eqref{eq:laplace_modbc}. \\

 \item With $u^{p}$ known on the boundary, we can solve \eqref{eq:doublelayerMu_disc} numerically with GMRES for $\mu$ at the Gauss--Legendre points, as described in subsection \ref{ss:homogeneuous}. The required modifications when $\Omega$ is a multiply connected domain are given in  \cite{Greenbaum1993}. Obtaining $u^{H}$ in $\Omega$ is just a matter of post processing, where the special quadrature is applied to points in $\Omega$ close to $\partial \Omega$. The final solution to \eqref{eq:poisson}--\eqref{eq:poissonBC} on a uniform grid in $\Omega$ is $u = u^{H} + u^{P}$. Note that we are not restricted to evaluating $u$ at these locations. Since the Fourier coefficients $\hat{u}^{P}$ and the density $\mu$ are known we can evaluate $u$ at any point in $\Omega$.
 
 \end{itemize}

\section{Numerical Results}
\label{s:numericalResults}
This section is organised as follows. We start with a general discussion of how to set the various parameters. Then  a strategy for finding appropriate values is presented, based on experiments for a simple numerical setting. Thereafter this strategy is shown to work also  for choosing parameters  for more advanced settings. Finally we present some timings to give an idea of the complexity of function extension in relation to solving the free--space Poisson equation.\\

 With the \emph{error} we refer to the relative discrete $\ell_{2}$ error for solving the full problem \eqref{eq:poisson}--\eqref{eq:poissonBC}, measured as
\begin{equation}
\frac{\|\mathbf{u}_{\text{exact}}-\mathbf{u}_{\text{numerical}}\|_{\ell_{2}}}{\|\mathbf{u}_{\text{exact}}\|_{\ell_{2}}},
\label{eq:relativeError}
\end{equation}
where
\begin{equation}
\|\mathbf{u}\|_{\ell_{2}} = \sqrt{\sum\limits_{i=1}^{N_{eval}^{2}}|u_{i}|^{2}/N_{eval}^{2}}
\label{eq:ell2}
\end{equation}
for a vector of length $N_{eval}^{2}$. We measure the error on an evaluation grid, which is a problem dependent uniform grid with resolution $2L/N_{eval}$, where $N_{eval}=1000$.\\

The domains used for numerical experiments are expressed in the complex plane as 

\begin{equation}
Re^{(n\, i\theta)} \left(c_{0} + \sum\limits_{j}(c_{j}\cos{(j\theta)}+d_{j}\sin{(j\theta)})\right) + a + ib
\label{eq:boundary}
\end{equation}
where $\theta\in[0,2\pi)$ and $n$ gives the orientation. The non--zero coefficients are stated for each complex multiply connected domain.\\

To obtain a numerical solution to the full problem \eqref{eq:poisson}--\eqref{eq:poissonBC}  the following parameters need to be set:
\begin{itemize}
\item $\varepsilon$: the shape parameter for the Gaussians used as basis for interpolation.
\item $R_{p}$: partition radius.
\item $N_{u}$: number of uniform grid points in one spatial direction.
\item $M$: The number of Vogel--nodes \eqref{eq:Vogelnodes}, i.e. the number of RBF centres per extension partition.
\item $N_{\partial \Omega}$: Number of Gauss--Legendre panels. Each panels has $16$ Gauss--Legendre points.

\item $L$: Length of the side of the box shaped computational domain $B$.
\item $\tilde{k}$: The regularity of the compactly supported RBF used to compute the weights \eqref{eq:weightfunction}.

\end{itemize}
The amount of overlap between partitions also has to be set. We always let them overlap by slightly  more than a radius $R_{p}$.  A large overlap gives the extension $f^{e}$  a more uniform band around $\Omega$ where it decays to zero, which grants a more predictable behaviour. Thus this parameter does not vary for different numerical settings. \\

In general, the solution is not sensitive to  the shape parameter $\varepsilon$, as long as $\varepsilon$ is within an appropriate range of values. We use $\varepsilon = 2$ for all numerical experiments in this article. Still, this is a problem dependent variable and may require tuning for substantially different domains and functions than presented in this paper. Observe that the right hand side \eqref{eq:results1} for example $1$ is an element of the native space $\mathcal{N}_{\varphi}$ for $\varepsilon = 2$, which is not the case for the right hand sides \eqref{eq:results2} and \eqref{eq:results3}, correspodning to example $2$ and $3$. \\

Another problem dependent parameter is the partition radius $R_{p}$. It is  strongly related to the grid resolution $N_{u}$ and the number of RBF centres $M$, i.e. the number of number of Vogel nodes \eqref{eq:Vogelnodes}. Recall that the RBF--QR algorithm has limited capacity for $M$, which is around $400$. Thus the greater the variation of $f$ the smaller partitions are required to obtain good approximations for the local interpolants. The choice of $M$ depends on the given $f$, and will thus be different for each numerical experiment.\\

 There is a simple way to find an appropriate value for $R_{p}$
 and to obtain a range for good choices of $\varepsilon$: measure the residual for solving the least squares problem \eqref{eq:LS} on some partition. A large residual implies that either $\varepsilon$ or $R_{p}$ needs to be changed. Another indication that smaller partitions are needed is if  $\max_{\mathbf{y}\in \mathbf{y}_{\Omega_{i},E}}{|s_{f,\Omega_{i}}(\mathbf{y})|}$ from \eqref{eq:extrapolatingf} is significantly larger than $\max_{\mathbf{y}\in \mathbf{y}_{\Omega_{i},\,\Omega}}{|f(\mathbf{y})|}$ for any partition $\Omega_{i}$; the extension should take on values in the same range as $f$.\\

  By construction $f^{e}$ decays to zero over a distance of  approximately $R_{p}$, meaning smaller partitions may require a larger value for $N_{u}$ to resolve $f^{e}$. Thus one may wish to set $R_{p}$ as large as possible, in order to avoid constructing an $f^{e}$ much harder to resolve than $f$. However, in most events $N_{u}$ is chosen to resolve $f^{e}$ close to the boundary of its support, not $f$ on $\Omega$. For difficult right hand sides, such as \eqref{eq:results2} from example $2$ and \eqref{eq:results3} from example $3$, the scale  is about the same. For simpler functions resolving $f^{e}$ is clearly harder, see for example \eqref{eq:results1}.  \\

  Previously we have said that \eqref{eq:LS} needs to be overdetermined, but no quantitative measure has been given of how much the solution is improved if more data points $\mathbf{y}_{\Omega_{i}}$, relative to $M$, is used. Therefore we introduce the measure
\begin{equation}
\label{eq:beta}
\beta_{i} = \frac{\text{number of data points in partition }\Omega_{i}}{M},
\end{equation}
and $\beta_{\min} = \min_{i}\beta_{i}$. It is an useful tool for analysing the relationship, in terms of error, between number of RBF centres $M$ and the available data in an extension partition, which is based on $N_{u}$.\\

The larger $\tilde{k}$, the higher regularity of $f^{e}$ and thus faster decay of the Fourier coefficients. But we also need to take into account that a larger $\tilde{k}$ implies that a finer grid is required to resolve the compactly supported RBF. In other words $\tilde{k}$ is strongly related to $R_{p}$ and $N_{u}$, just as $M$ is. In the following subsection we devise a scheme via an heuristic approach to choose  $\tilde{k}$ and $M$. It is based on the measure of the \emph{number of uniform grid points per partition radius}, which we denote $P$ and is given by
\begin{equation}
P = \frac{N_{u}}{2 L} R_{p}.
\label{eq:pntsperRp}
\end{equation}

The length $L$ of the sides of the box $B$ is set such that all
$\mathbf{y}\in\cup_{i}\bar{\Omega}_{i}$ are in $B$. The number of Gauss--Legendre panels $N_{\partial \Omega}$ is chosen  large enough to resolve $g$, which in most events is sufficient to resolve $g-u^{p}$ as well. This is  verified numerically in subsections \ref{ss:ex2} and \ref{ss:ex3}.

\subsection{Example $1$: Parameter selection}
\label{ss:parameterSelection}

As right hand side for \eqref{eq:poisson}--\eqref{eq:poissonBC}, consider the smooth function
\begin{equation}
\label{eq:results1}
f(x,y) = -\sin{(2\pi x)}\sin{(2\pi y)},
\end{equation}
defined in a disc centred at ${(17/701,\,5/439)}$ with radius one. For this simple example we will study different errors and discuss how to choose parameters.  \\

\begin{figure}[!ht]
\centering\includegraphics[width=0.5\linewidth]{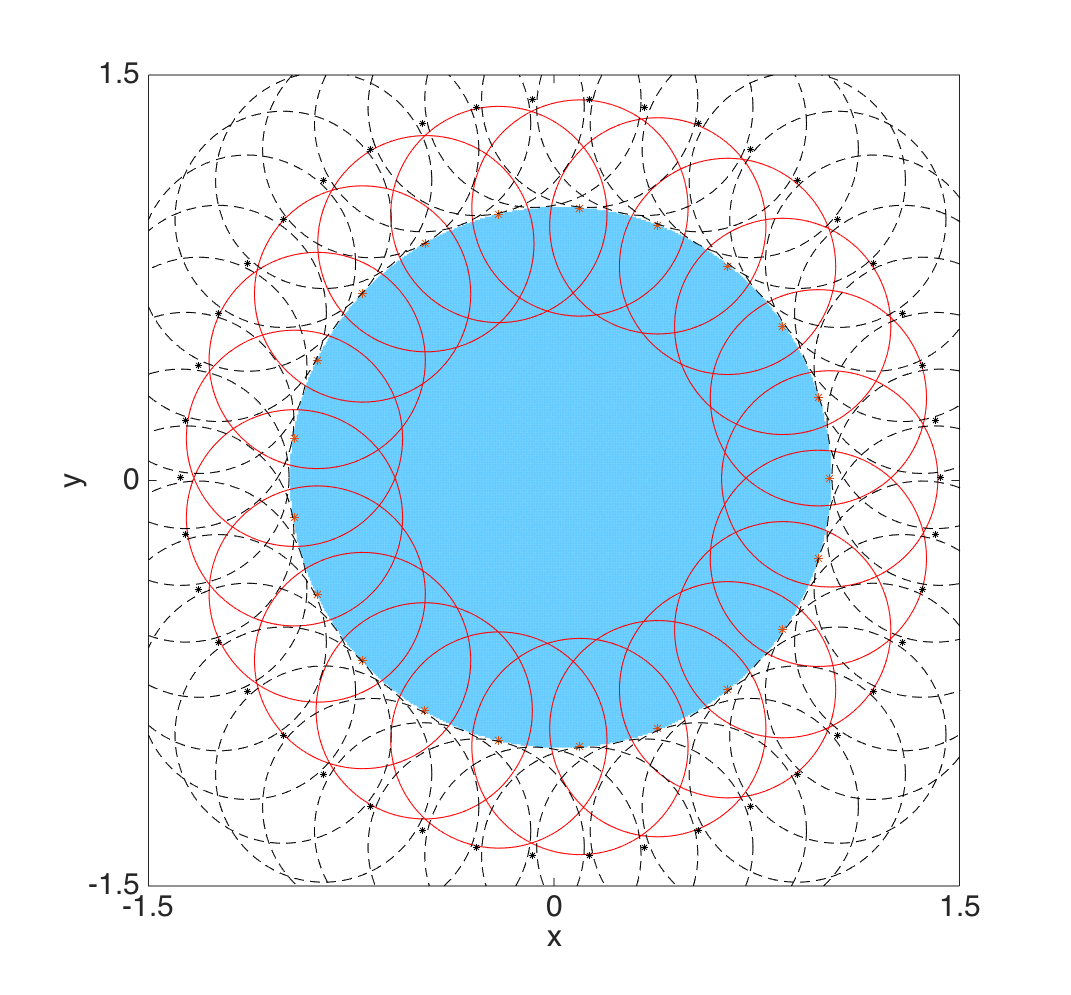}
\caption{Regular covering of the disc $\Omega$ used in subsection \ref{ss:parameterSelection}. Red circles are extension partitions, black are zero--partitions. Red and black stars are centres for the extension partitions and zero partitions, respectively. }
\label{fig:results1covering}
\end{figure}

Following the approach given above we set: $L = 1.5$ and the partition radius $R_{p}$ to  $0.4$ . Thus we use $21$ extension partitions and $42$ zero--partitions distributed as in Figure \ref{fig:results1covering}. The boundary is discretised with $32$ Gauss--Legendre panels. The resolution $N_{u}$ varies from $40$ to $500$, and for each value the parameters $M$ and $\tilde{k}$ need to be set.\\

\FloatBarrier


We begin by investigating the influence of the regularity of the weight function on the error convergence. To isolate this error, we want to remove the error from the local extension. Hence, instead creating an extension of $f$ as in \eqref{eq:extrapolatingf} we pick a smooth $f$ defined in all of $B$ and set for each interpolation patch $s_{f,\Omega_{i}} = f$. Only the choice of weight function, i.e. $\tilde{k}$, and the resolution of the uniform grid will change, the rest will be fixed. \\
\begin{figure}[ht]
    \begin{center}
        \includegraphics[width=0.49\textwidth,height=7cm]{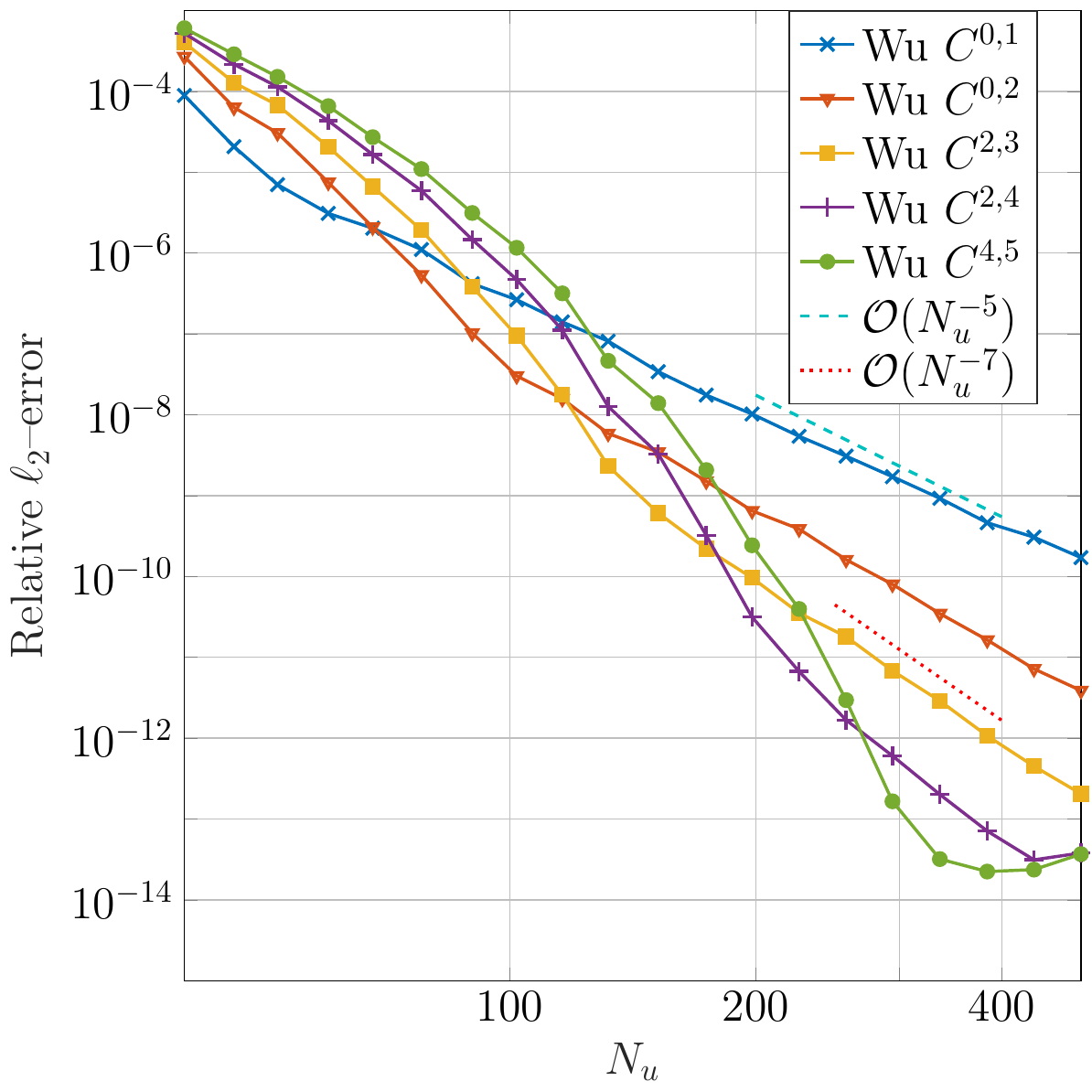}
        \includegraphics[width =0.49\textwidth,height=7cm]{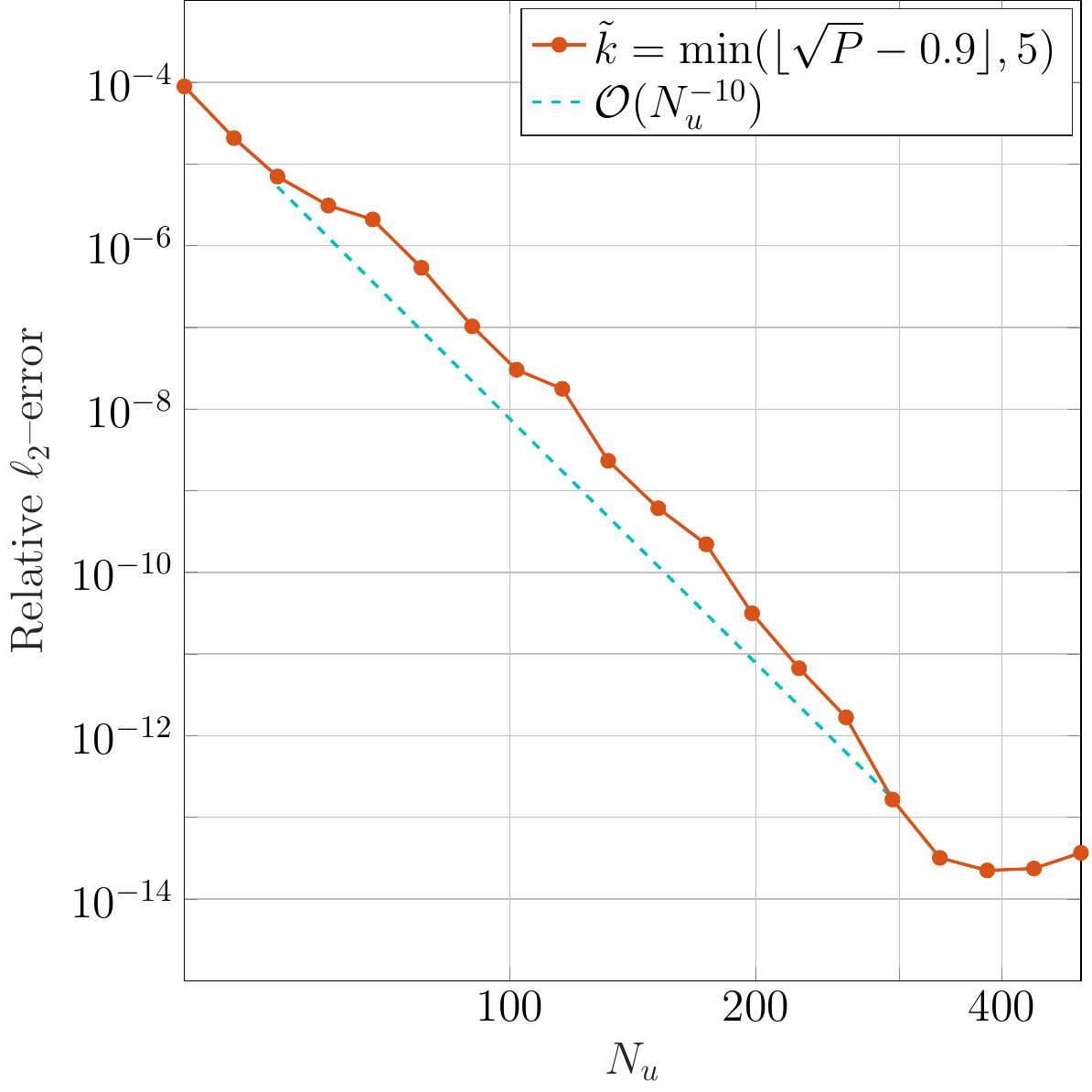}
    \end{center}
    \caption{Error in numerical solution for the Poisson equation with right hand side given by \eqref{eq:results1},  but with local extensions given by analytic expression. Left: relative $\ell_{2}$ error as a function of $N_{u}$ in loglog-scale for various compactly supported RBF given in Table \ref{tab:Wu}. Right: relative $\ell_{2}$ error as a function of $N_{u}$ in loglog-scale with $\tilde{k}$ set by \eqref{eq:ktilde}.}
 \label{fig:results1simpleAndWu}
\end{figure}
\begin{figure}[ht]
    \begin{center}
        \includegraphics[width=0.49\textwidth,height=7cm]{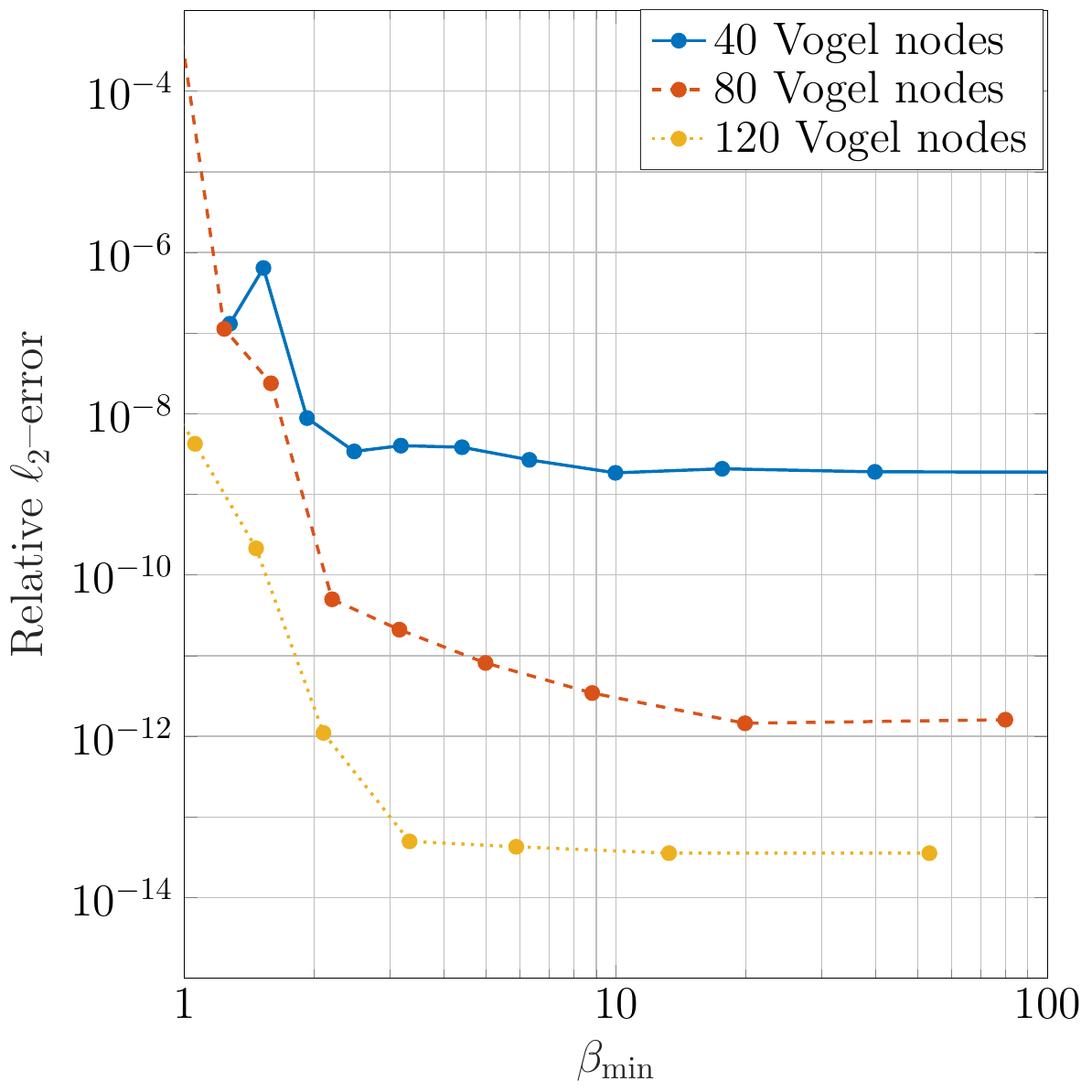}
        \includegraphics[width =0.49\textwidth,height=7cm]{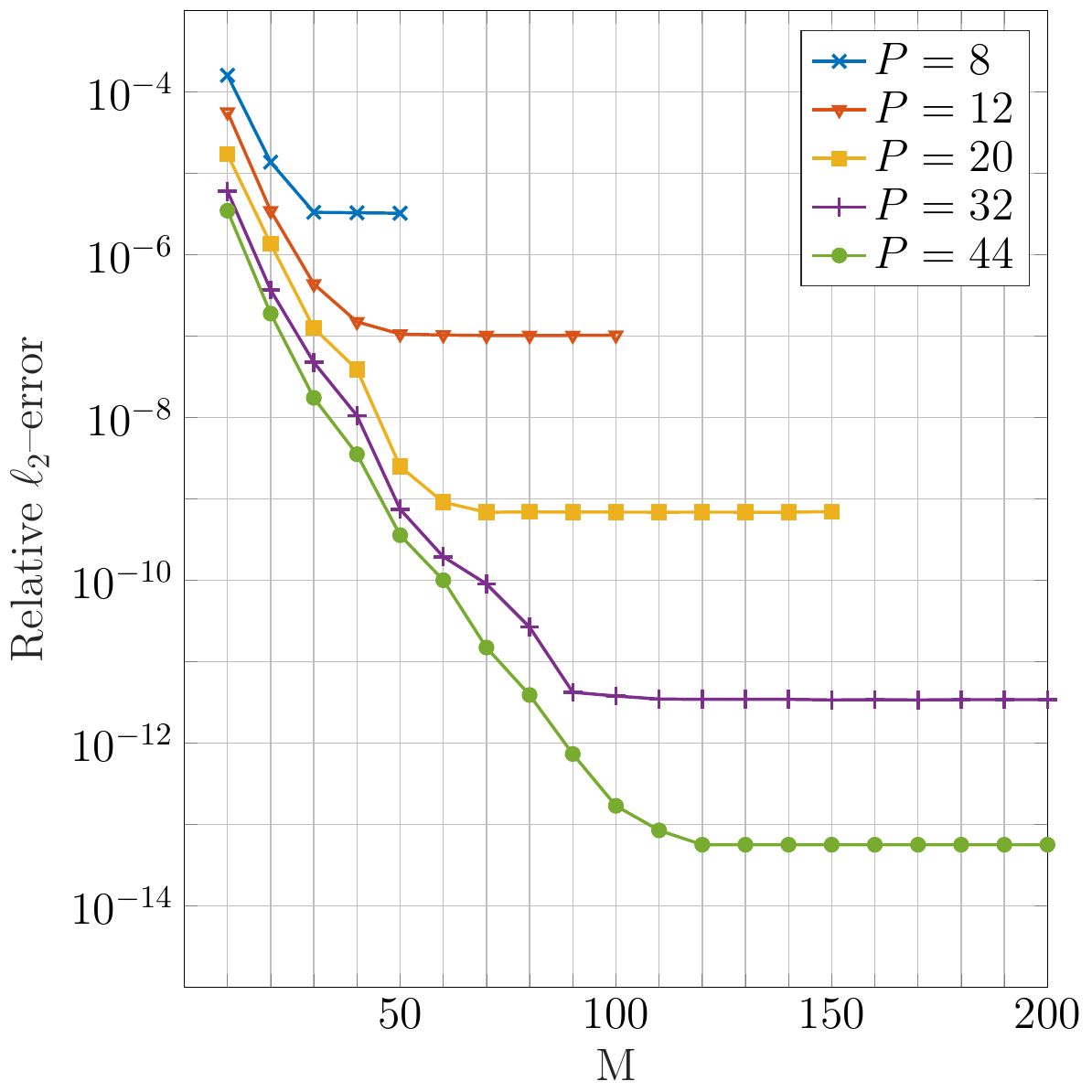}
    \end{center}
    \caption{Error in numerical solution for the Poisson equation  with $f$ given by \eqref{eq:results1}. Effect of parameter choices for computing local extension on each patch, with $\tilde{k}$ set by \eqref{eq:ktilde}. Left: Relative $\ell_{2}$ error as a function of $\beta_{\min}$, defined below \eqref{eq:beta}, in semilog-scale. Right: Relative $\ell_{2}$ error as function of number of Vogel nodes $M$, for different values of $P$ \eqref{eq:pntsperRp}.}
 \label{fig:results1betaandM}
\end{figure}
\FloatBarrier
In the left plot of Figure \ref{fig:results1simpleAndWu} the convergence of the relative $\ell_{2}$ error is plotted in loglog--scale as a function of $N_{u}$ for various compactly supported RBF used to construct $w$. Recall that an RBF in $C^{k,\,\tilde{k}}$ is of regularity $k$ at origin and $\tilde{k}$ at the boundary of its support. Initially the error for $C^{2,3}$  converges  super algebraically and is followed by a tail with algebraic convergence, as expected. The order of convergence for the algebraic tail is inherited by the the RBFs regularity at the edge of its support. However, this is only true if we never evaluate any RBF in a neighbourhood around its origin. Compare the curves corresponding to Wu $C^{0,1}$ and Wu $C^{0,2}$ in the left plot in Figure \ref{fig:results1simpleAndWu}. The former has  an algebraic tail with slope $-5$, while the latter has $-6$. For Wu $C^{0,1}$ the error is dominated by the algebraic tail of order $5$ for the entire spectrum of $N_{u}$. This indicates that the error associated with resolving the weights $w$ is largest. We observe that the error for Wu function $C^{4,5}$  converges super algebraically down to a relative error of $10^{-14}$, but requires a fine grid resolution to do so. First at approximately $N_{u} = 260$ is Wu function $C^{4,5}$ the better choice; this corresponds to $P \approx 35$ . These results suggest that a higher rate of convergence can be obtained by picking an optimal RBF to construct the weight function for each $N_{u}$, compared to using the same for all $N_{u}$.  Heuristically we have found that with
\begin{equation}
\label{eq:ktilde}
\tilde{k} = \min{\left(\Bigl\lfloor\!\sqrt{P}-0.9\Bigl\rfloor,5\right)},
\end{equation}
where $\lfloor x\rfloor$ gives the greatest integer less than or equal to $x$, we essentially obtain algebraic convergence corresponding to $\mathcal{O}(N_{u}^{-10})$. This is evident from the right plot in Figure \ref{fig:results1simpleAndWu}.    Various numerical simulations confirm that the optimal choice for a given $N_{u}$ is roughly the same even for different length scales and other functions $f$. This will be shown in subsections \ref{ss:ex2} and \ref{ss:ex3}. Unless stated otherwise, in subsequent numerical experiments $\tilde{k}$ is given by \eqref{eq:ktilde}.\\

From now on we the create local extension by \eqref{eq:functionextension}, meaning $M$ needs to be set. The left plot in Figure \ref{fig:results1betaandM} shows the error \eqref{eq:relativeError} as function of $\beta_{min}$. Typically $\beta_{min}=3$ is  sufficient, and little is gained by increasing the number of data locations further. Given a fine uniform grid one can thus downsample for the interpolation problem. Solving a least squares problem is executed efficiently in \textsc{Matlab}, for reference: computing $A\backslash Y$ for $A\in\mathbb{R}^{600\times 200}$, $Y\in\mathbb{R}^{600}$ with randomised elements takes approximately $0.0054$ seconds on a MacBook Pro with a $2.6$ GHz Intel Core i$5$ processor and $16$ GB of $1600$ MHz DDR$3$ memory. \\

In Figure \ref{fig:results1betaandM} the relative $\ell_{2}$ error is plotted as a function of $M$ for different $P$. Clearly the error decreases exponentially initially, but then level out around $M\approx 3P$. As safety measure we use $4P$, but note that for small $N_{u}$ the resulting $\beta_{\min}$ may be less than three. In such case, we adjust $M$ accordingly: the number of uniform grid points inside $\Omega$ in a given partition can be approximated by $\pi P^{2}/2$. To obtain a least--squares system overdetermined by at least a factor of two, $M$ should be less than half of that, e.g.  $(0.8/2)\pi P^{2}/2$. Thus the scheme we apply is
\begin{equation}
\label{eq:autoM}
M = \min{(0.8\pi P^{2}/4,4P)},
\end{equation}
although a majority of numerical settings yield $M = 4P$, since $0.8\pi P^{2}/4 = 4P$ for $P\approx 6.4$. Note that  neither \eqref{eq:ktilde} nor \eqref{eq:autoM} are optimal, but meant as a guide to select the parameters.\\

\begin{figure}[ht]
    \begin{center}
        \includegraphics[width=0.49\textwidth,height=7cm]{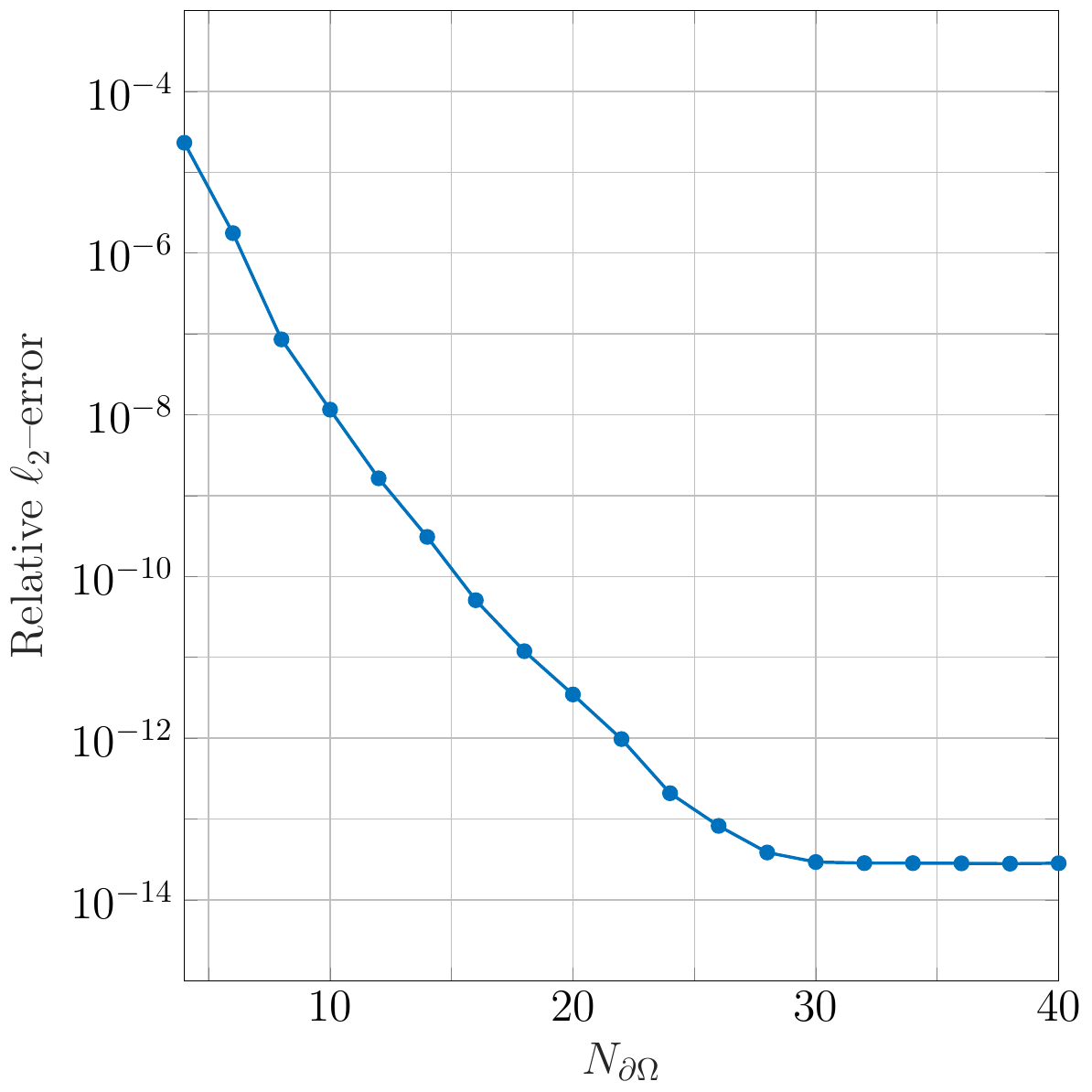}

    \end{center}
    \caption{Relative $\ell_{2}$ error as a function of the number of Gauss--Legendre panels $N_{\partial \Omega}$ in semilog-scale. Results are for solving the Poisson equation with $f$ given by \eqref{eq:results1}, with local extensions given by analytic expression.}
 \label{fig:results1GL}
\end{figure}

\FloatBarrier


In Figure \ref{fig:results1GL} the relative $\ell_{2}$ error is plotted against the number of Gauss--Legendre panels. Thus we let $N_{\partial\Omega}$ vary for $N_{u} = 400$ and other parameters chosen as before. Since the amount of panels sets the discrete representation for the boundary, which is used to sort points as in $\Omega$ or $E$, the evaluation grid is pruned. Only points considered in $\Omega$ for the all the investigated resolutions are kept. We see the expected rapid convergence as the number of panels is increased.\\



\subsection{Example $2$: The Poisson equation on a multiply connected domain}
\label{ss:ex2}

 We replicate the setting for the most difficult numerical experiment performed in \cite{ASKHAM20171}, but give the details here for convenience. For this experiment we study only the error \eqref{eq:relativeError} as a function of $N_{u}$. We consider the right hand side 
\begin{equation}
\label{eq:results2}
f(x,y) = -200\sin{(10(x+y))} + \frac{2}{9} + 1000(1000x^2-1)e^{-500x^2}
\end{equation}
for the Poisson equation \eqref{eq:poisson}--\eqref{eq:poissonBC} on a complex multiply connected domain: for the outer boundary the non--zero coefficients for \eqref{eq:boundary} are $c_{0} = 0.25$, $d_3 = c_{6} = c_{8} = c_{10} = 0.01$, $c_{5} = 0.02$, $R=1$and $n = 1$. For the inner boundary we set $c_{0} = 0.05$, $c_{2} = d_{3} = c_{5} = c_{7} = 0.005$, $R=1$ and $n=-1$.\\

 \begin{figure}[ht]
    \begin{center}
        \includegraphics[width=0.49\textwidth]{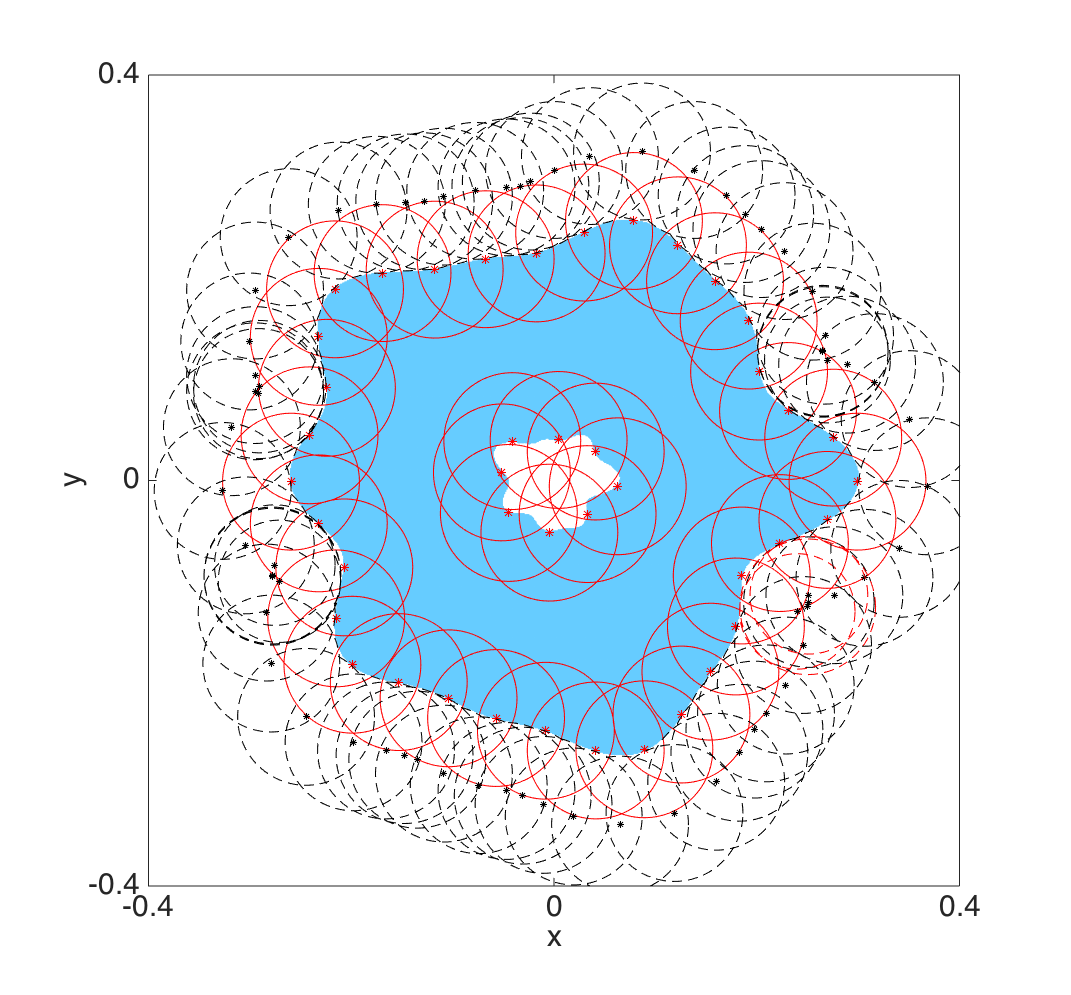}
        \includegraphics[width =0.49\textwidth]{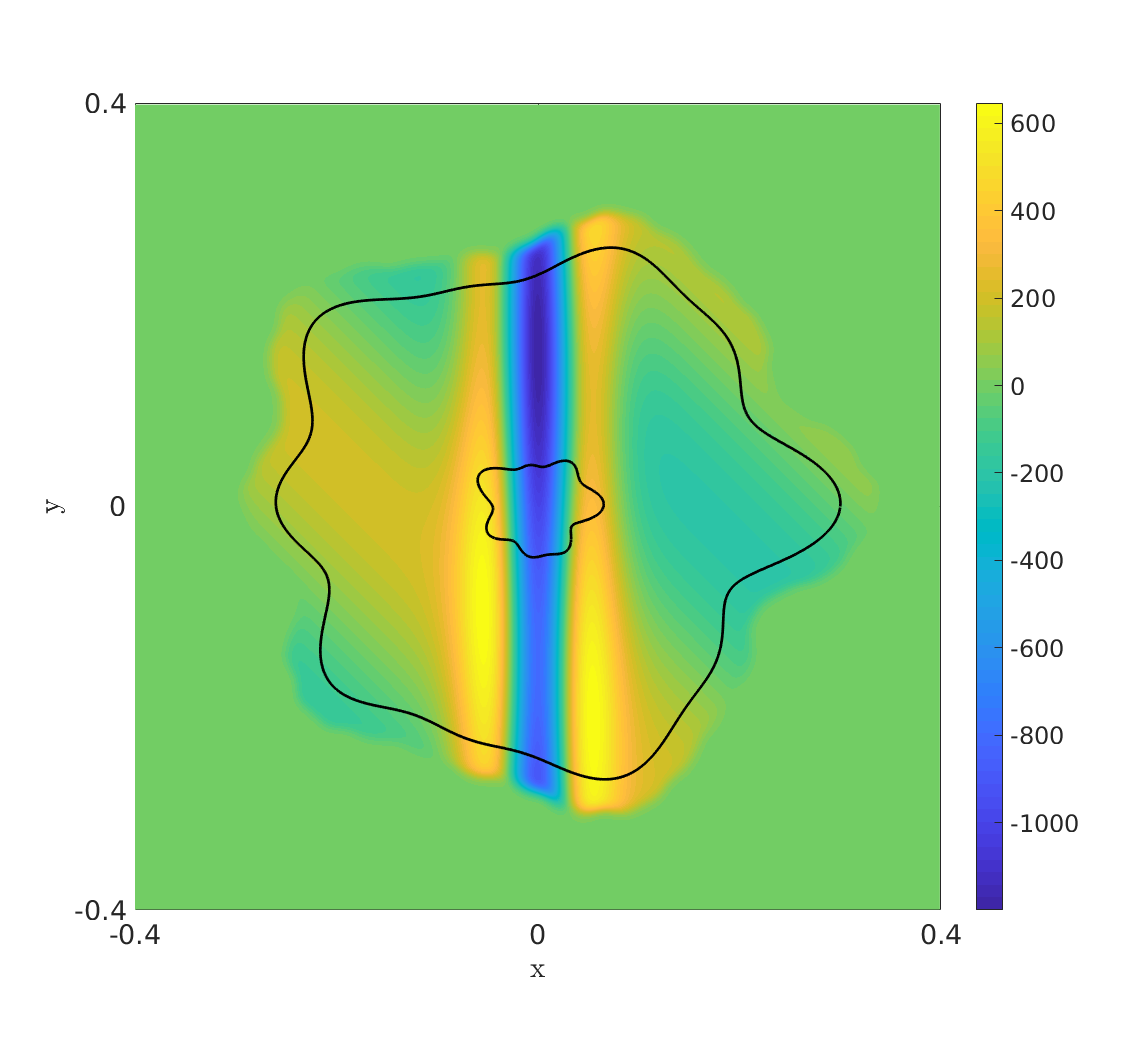}
    \end{center}
    \caption{Left: Regular covering of $\Omega$ from subsection \ref{ss:ex2}. Red circles are extension partitions, black are zero--partitions. Red and purple black are centres for the extension partitions and zero partitions, respectively. Right: Extension $f^{e}$ by PUX of $f$, given by \eqref{eq:results2}. }
    \label{fig:example2Coverfe}
\end{figure}

 The given right-hand side \eqref{eq:results2} exhibits steep edges parallel to the $y$--axis and is visualised in Figure \ref{fig:example2Coverfe}. The covering consists of $38$ non--zero partitions on the outer boundary and $9$ partitions on the inner boundary. Different radius for the extension partitions along the outer boundary and inner boundary can be used, which requires the precomputation  two different matrices $A$. For simplicity we let all extension matrices have the same radius, therefore they overlap more along the boundary of the cavity. The regular covering of $\Omega$ is seen in Figure \ref{fig:example2Coverfe}.  \\

\begin{figure}[ht]
    \begin{center}
        \includegraphics[width=0.49\textwidth]{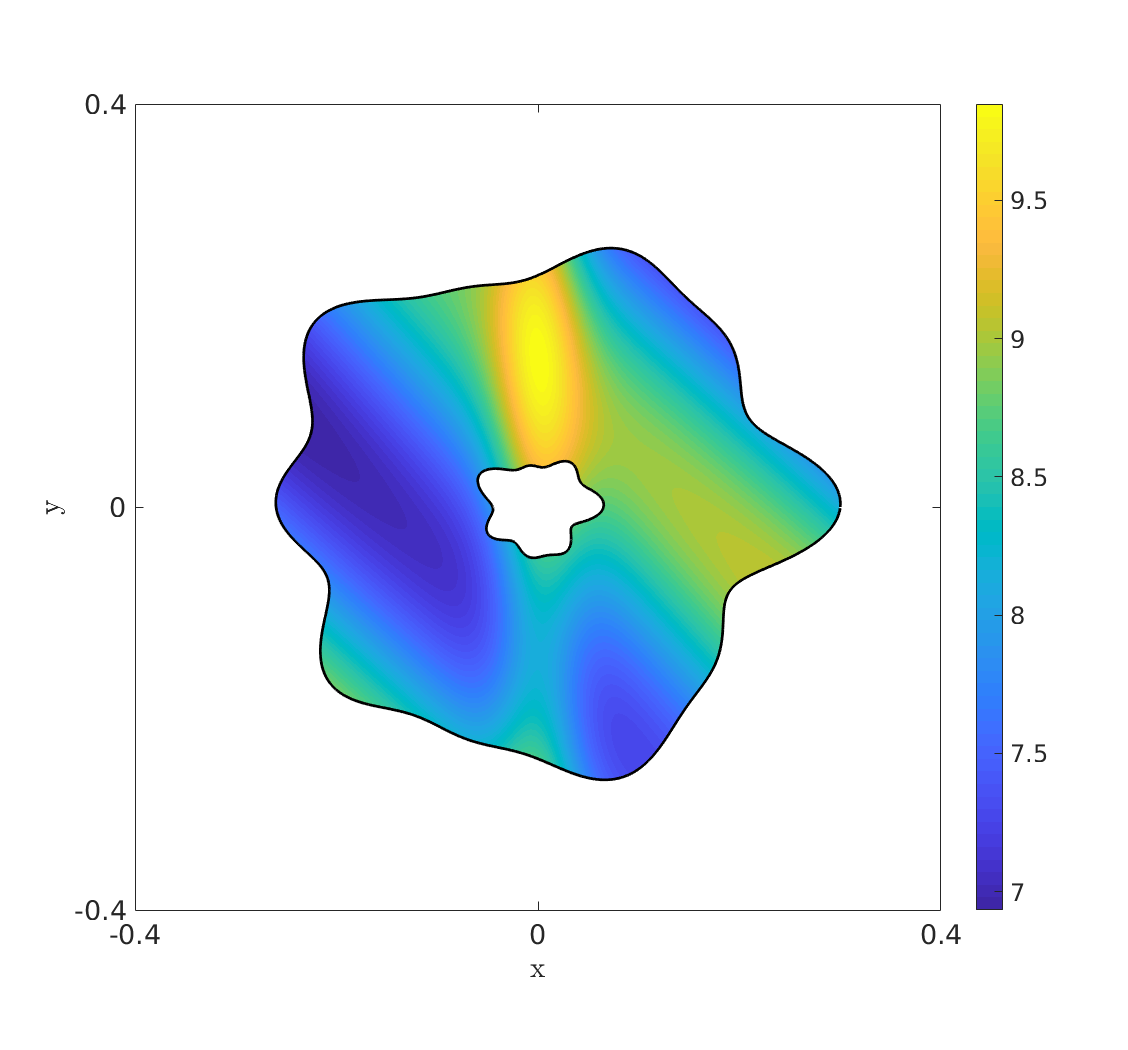}
        \includegraphics[width =0.49\textwidth]{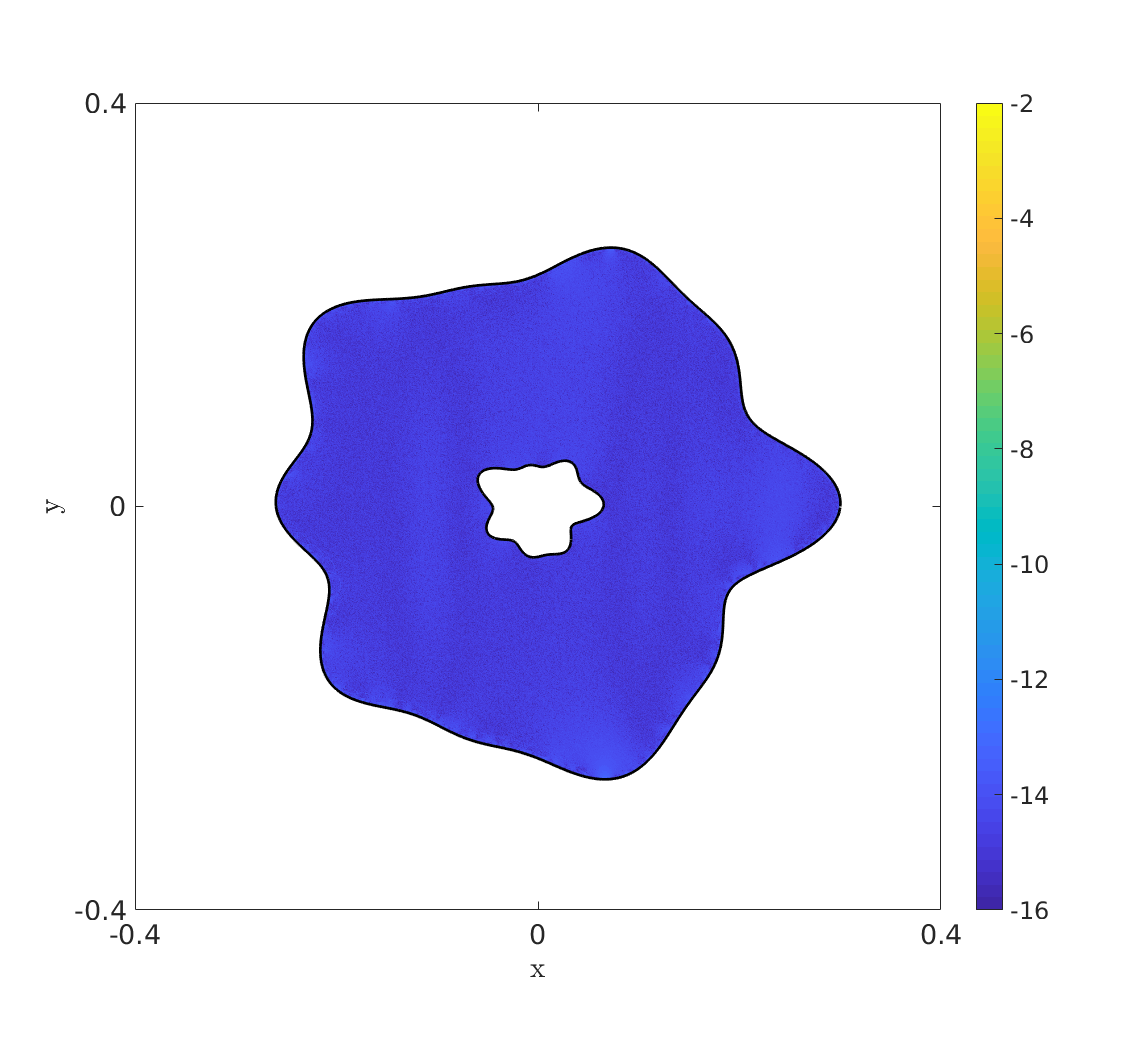}
    \end{center}
    \caption{ Left: Numerical solution $u$ to the Poisson  equation with $f$ from \eqref{eq:results2} for $N_{u} = 700$. Right: Pointwise relative $\log_{10}$--error for $N_{u} = 700$. Max relative error is $3.06\times 10^{-13}$. }
    \label{fig:example2uanderror}
\end{figure}

For the numerical setting we let $L = 0.4$, use $64$ and $44$ Gauss--Legendre panels for the outer boundary and for the inner boundary, respectively. The radius for the extension partitions $R_{p}$ is  $0.0675$. The Poisson equation is solved for grids with $N_{u}$ ranging from $10^2$ to $10^3$. We apply \eqref{eq:ktilde} and \eqref{eq:autoM} to pick $\tilde{k}$ and $M$. An extension $f^{e}$ of \eqref{eq:results2} can be seen in Figure \ref{fig:example2Coverfe}, corresponding to $N_{u} = 700$. The suppression to zero is rapid outside the exterior boundary of $\Omega$ and almost looks discontinuous. The extension of $f$ into the cavity seems to mimic $f$ remarkably well. However, a few remarks concerning the extension $f^{e}$ in the cavity of $\Omega$ are in order: It is not obvious that no zero partitions are required in the cavity to obtain a well--behaved extension. Recall that far away from the boundary the extension can grow dramatically. This is usually an indication that smaller partitions are needed. Here we ended up with a setting where every point in the cavity belongs to an extension partition and the extension is well--behaved. Thus there is no need for zero partitions, whereas they will indeed be needed in example $3$.\\

\begin{figure}[ht]
    \begin{center}
        \includegraphics[width=0.49\textwidth,height=7cm]{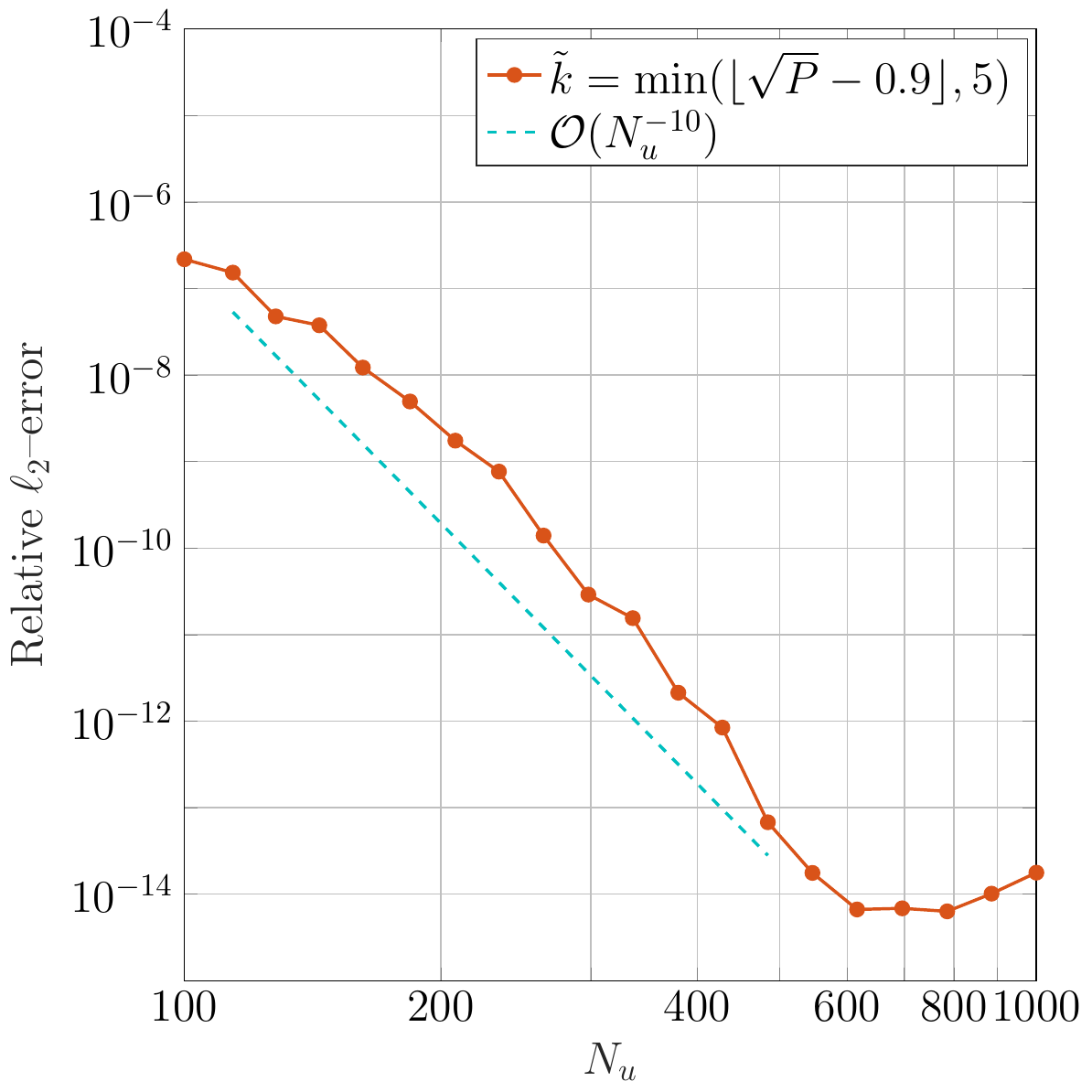}
        \includegraphics[width =0.49\textwidth,height=7cm]{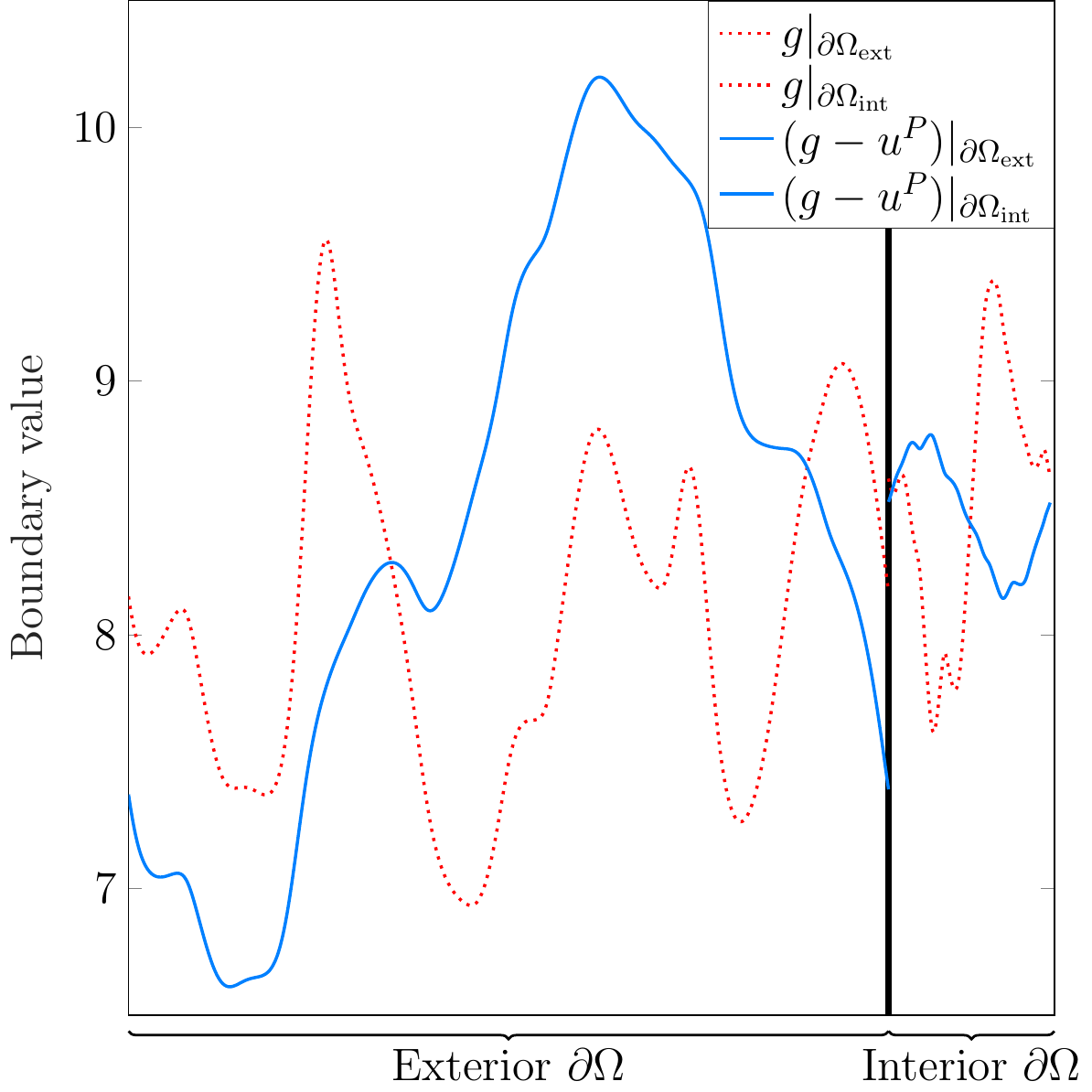}
    \end{center}
    \caption{Results for solving the Poisson equation with \eqref{eq:results2}. Left: Relative error as function of uniform grid resolution $N_{u}$, with $\tilde{k}$ and $M$ set by \eqref{eq:ktilde} and \eqref{eq:autoM}. Right: The boundary values $g$ and modified boundary values $g-u^{P}$ on the exterior and interior boundary segments.}
    \label{fig:function2convergenceAndbc}
\end{figure}

For the  resolution $N_{u} = 700$ the numerical solution $u$  and the corresponding pointwise relative error are plotted in Figure \ref{fig:example2uanderror}. In Figure \ref{fig:function2convergenceAndbc} a convergence plot as a function of $N_{u}$ is shown. As for the simpler setting we obtain an order of convergence corresponding to $\mathcal{O}(N_{u}^{-10})$. This is due to the shift to a Wu function of higher regularity at suitable points as $N_{u}$ increases. Thus the estimate \eqref{eq:ktilde} appears to be useful for more difficult numerical settings as well. For comparison, Askham et al. in \cite{ASKHAM20171} obtained a relative max error around $10^{-6}$ for an adaptive grid with in total $10^{6}$ points with, convergence of order $3$. In their embedded boundary-setting  a $C^{0}$--extension is constructed by solving the Laplace equation with $f$ as boundary condition, as briefly explained in the introduction. 
\\

In Figure \ref{fig:function2convergenceAndbc} the boundary condition $g$ and modified boundary conditions $g-u^{P}$ are plotted, where the latter is passed as input to the integral equation solver. We observed in several simulations that $g-u^{P}$ is not harder to resolve than $g$, meaning one can choose the discretisation of $\partial \Omega$ a priori based on the originally stated boundary conditions $g$.\\

\FloatBarrier

\subsection{Example $3$: The Poisson equation on a domain with a larger cavity and more oscillating right hand side}
\label{ss:ex3}

The third test features the right--hand side 
\begin{equation}
\label{eq:results3}
f(x,y) = -\sum_{i = 0}^{5} 2^{2i}e^{-\sqrt{2^{i}}}(\cos{(2^{i}x)} + \cos{(2^{i}y)})
\end{equation}
for the Poisson equation \eqref{eq:poisson}--\eqref{eq:poissonBC}. It oscillates with high frequency around the boundaries and is thus overall harder to resolve than \eqref{eq:results2}. The domain is again multiply connected and the non--zero coefficients for the outer boundary are $c_{0} = 1$, $c_{-5} = d_{-1} = 0.2$, $R=1$ and $n = 1$. For the inner boundary we set $c_{0} = 1$, $c_{-6} = d_{-3} =  0.1$, $R=0.3$, $b = 0.17$ and $n=-1$.  \\

\begin{figure}[ht]
    \begin{center}
        \includegraphics[width=0.49\textwidth]{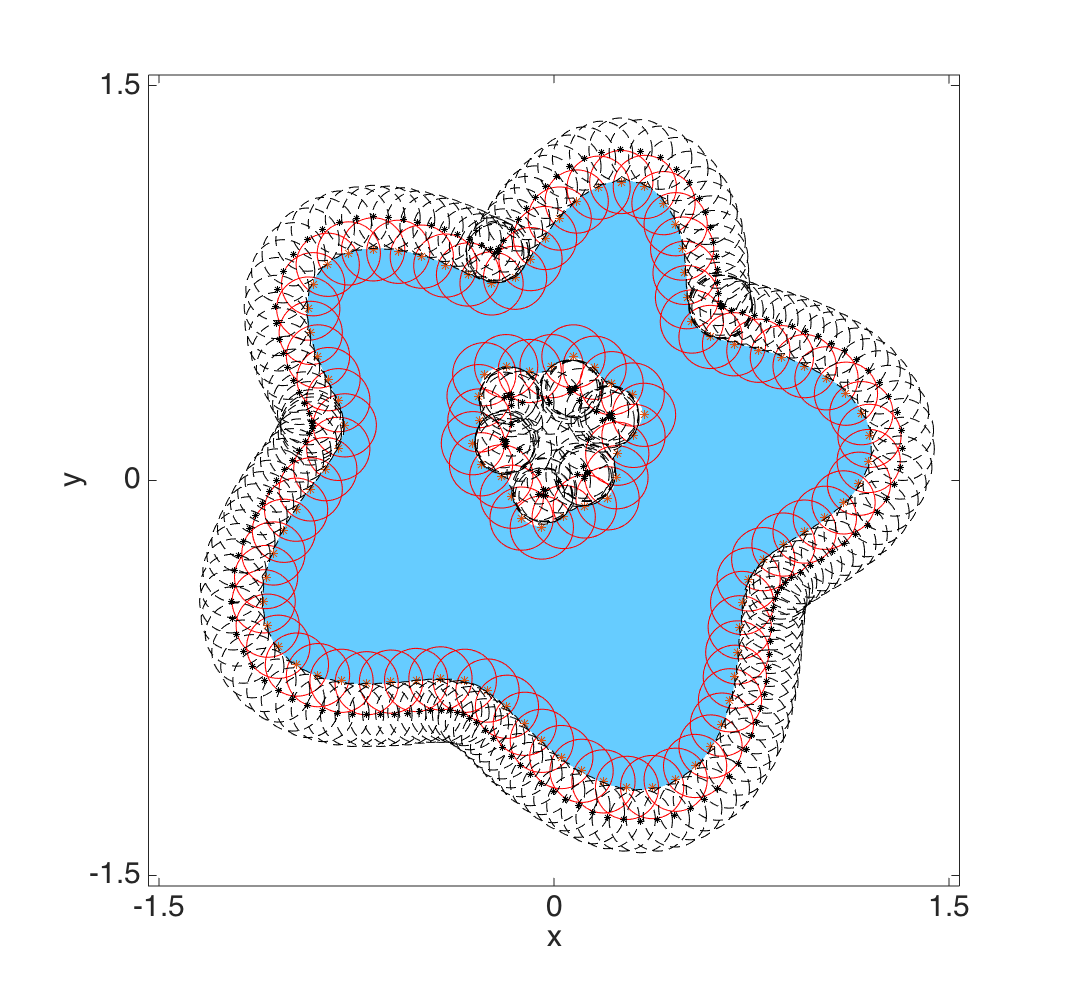}
        \includegraphics[width =0.49\textwidth]{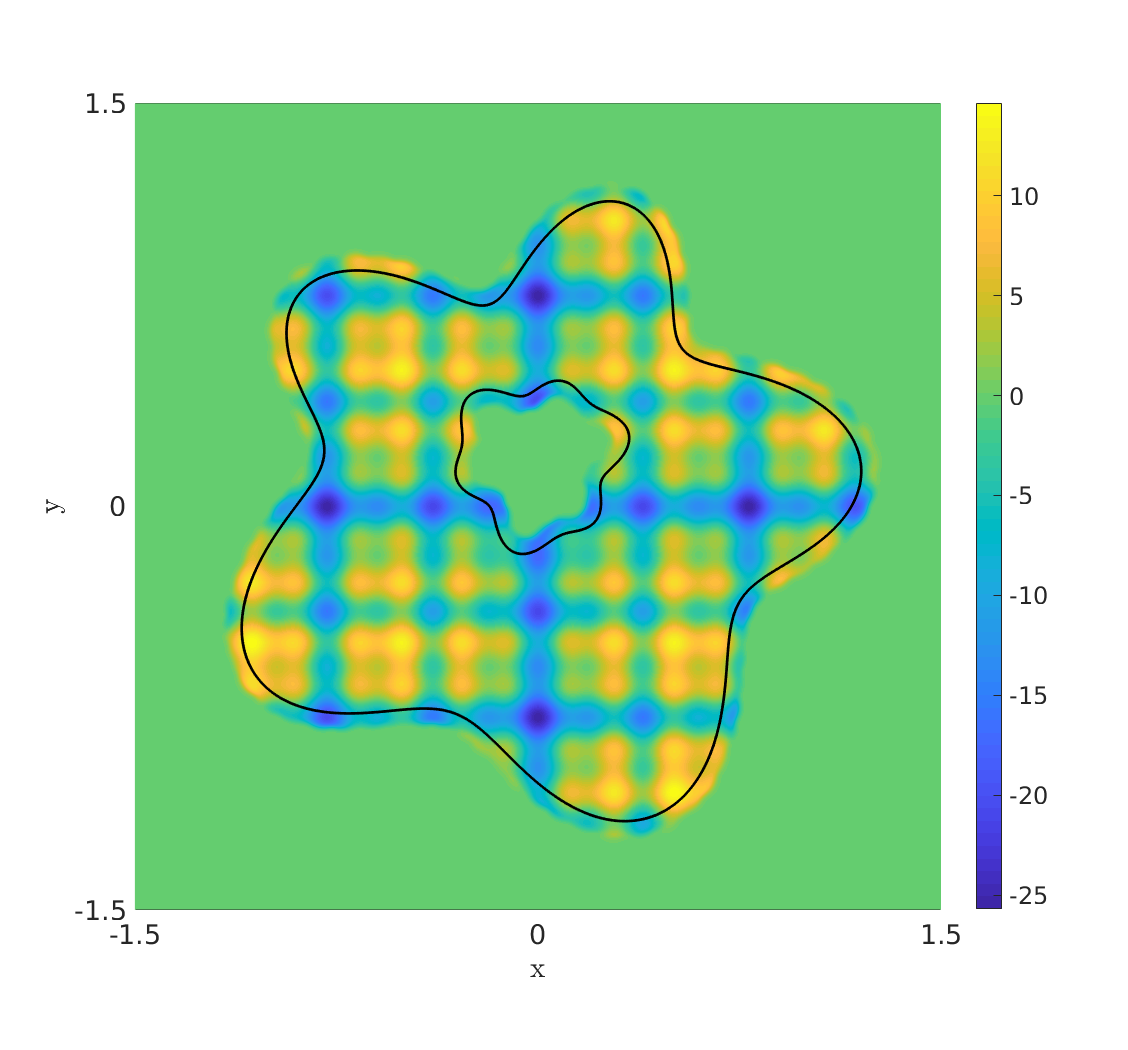}
    \end{center}
    \caption{Left: Regular covering of $\Omega$ from subsection \ref{ss:ex3}. Red circles are extension partitions, black are zero--partitions. Right: Extension $f^{e}$ by PUX of $f$, given by \eqref{eq:results3}. }
    \label{fig:example3Coverfe}
\end{figure}

We use $L = 1.54$, $84$  Gauss-Legendre panels and $82$ partitions for the outer boundary and $40$ Gauss-Legendre panels and $23$ partitions for the inner boundary.   Once again $N_{u}$ ranges from $10^2$ to $10^3$, $\tilde{k}$ and $M$ is set by \eqref{eq:ktilde} and \eqref{eq:autoM} and all extension partitions have the same radius, $R_{p} = 0.12$. Consequently the partitions along the inner boundary do not cover the cavity, thus $f^{e}$ must be suppressed to zero inside. See Figure \ref{fig:example3Coverfe} for $f^{e}$ and the covering, which should be compared with the covering from the previous example, see Figure \ref{fig:example2Coverfe}. \\

\begin{figure}[ht]
    \begin{center}
        \includegraphics[width=0.49\textwidth]{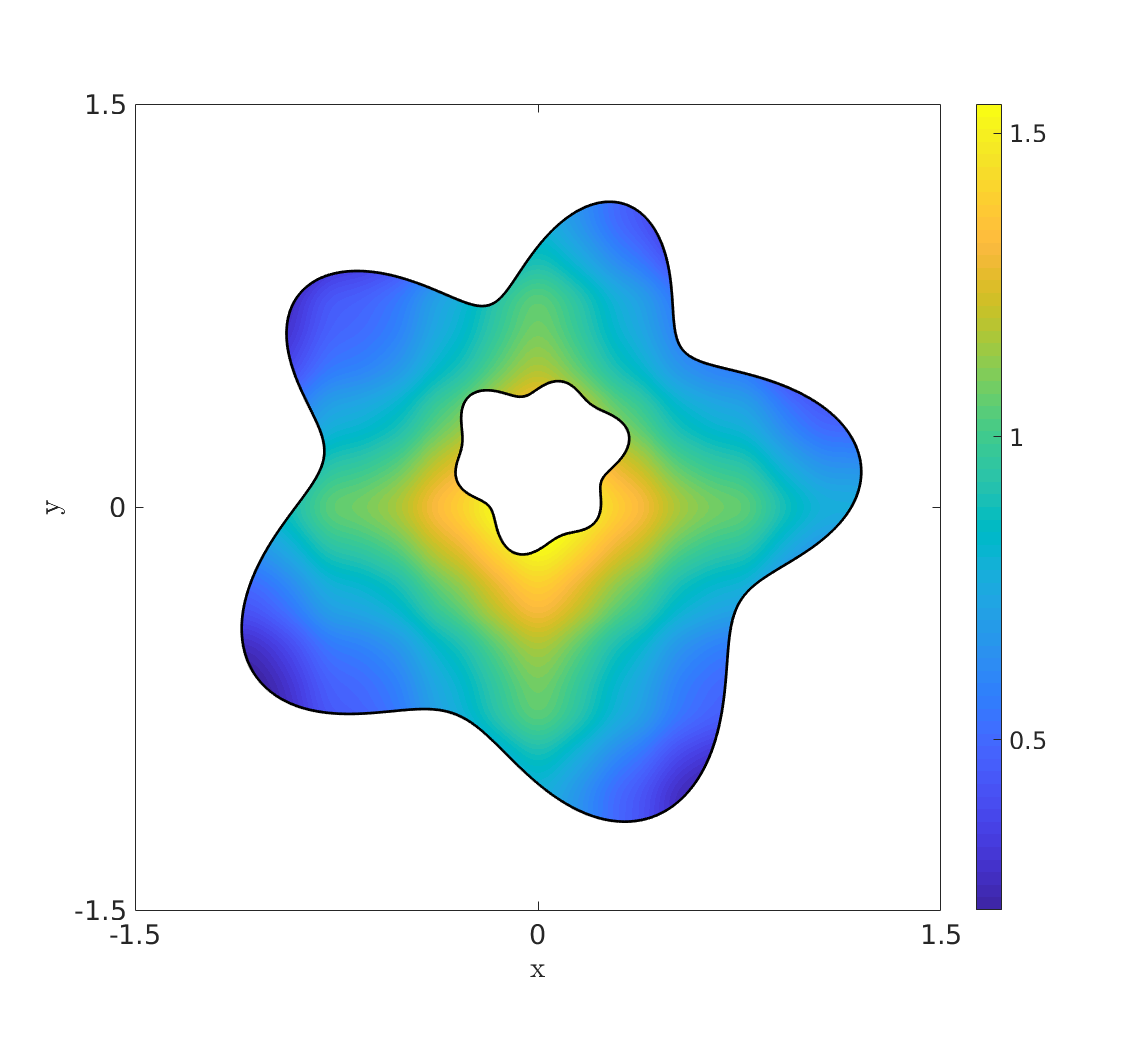}
        \includegraphics[width =0.49\textwidth]{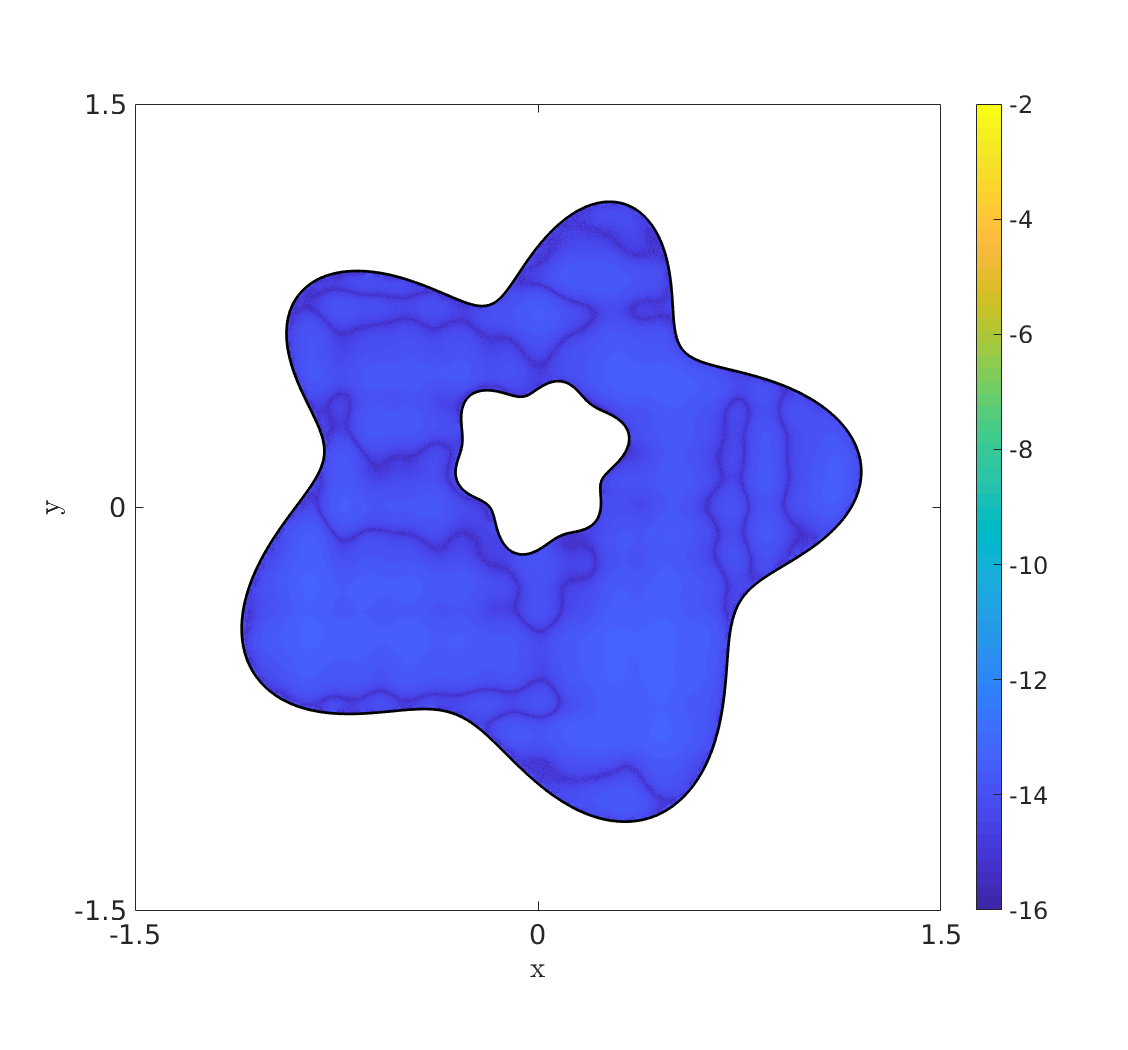}
    \end{center}
    \caption{ Left: Numerical solution $u$ to the Poisson equation with f from \eqref{eq:results3} for $N_{u} = 1000$. Right: Pointwise relative $\log_{10}$--error for $N_{u} = 1000$. Max relative error is $3.53\times 10^{-13}$.}
    \label{fig:example3uanderror}
\end{figure}

\begin{figure}[ht]
    \begin{center}
        \includegraphics[width=0.49\textwidth,height=7cm]{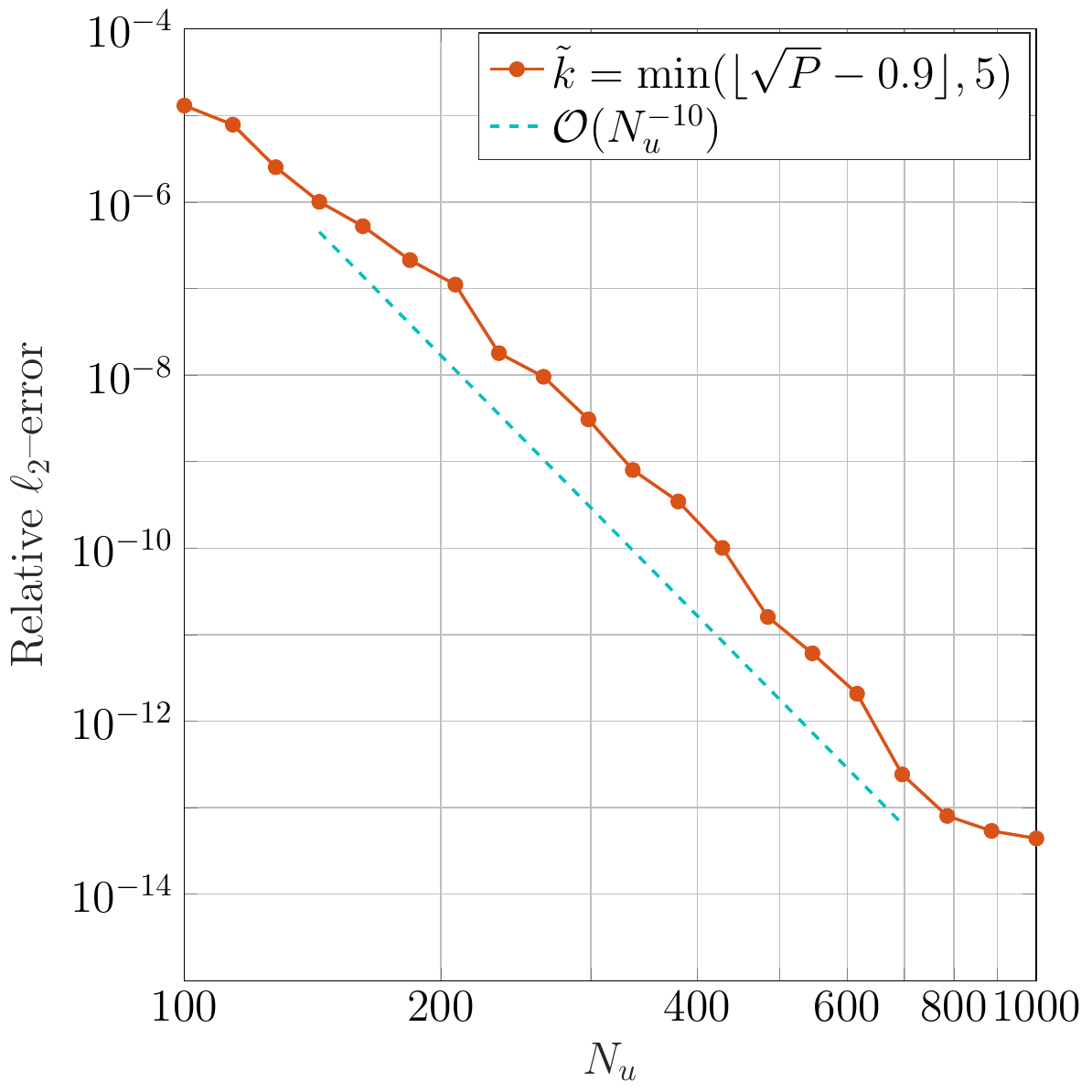}
        \includegraphics[width =0.49\textwidth,height=7cm]{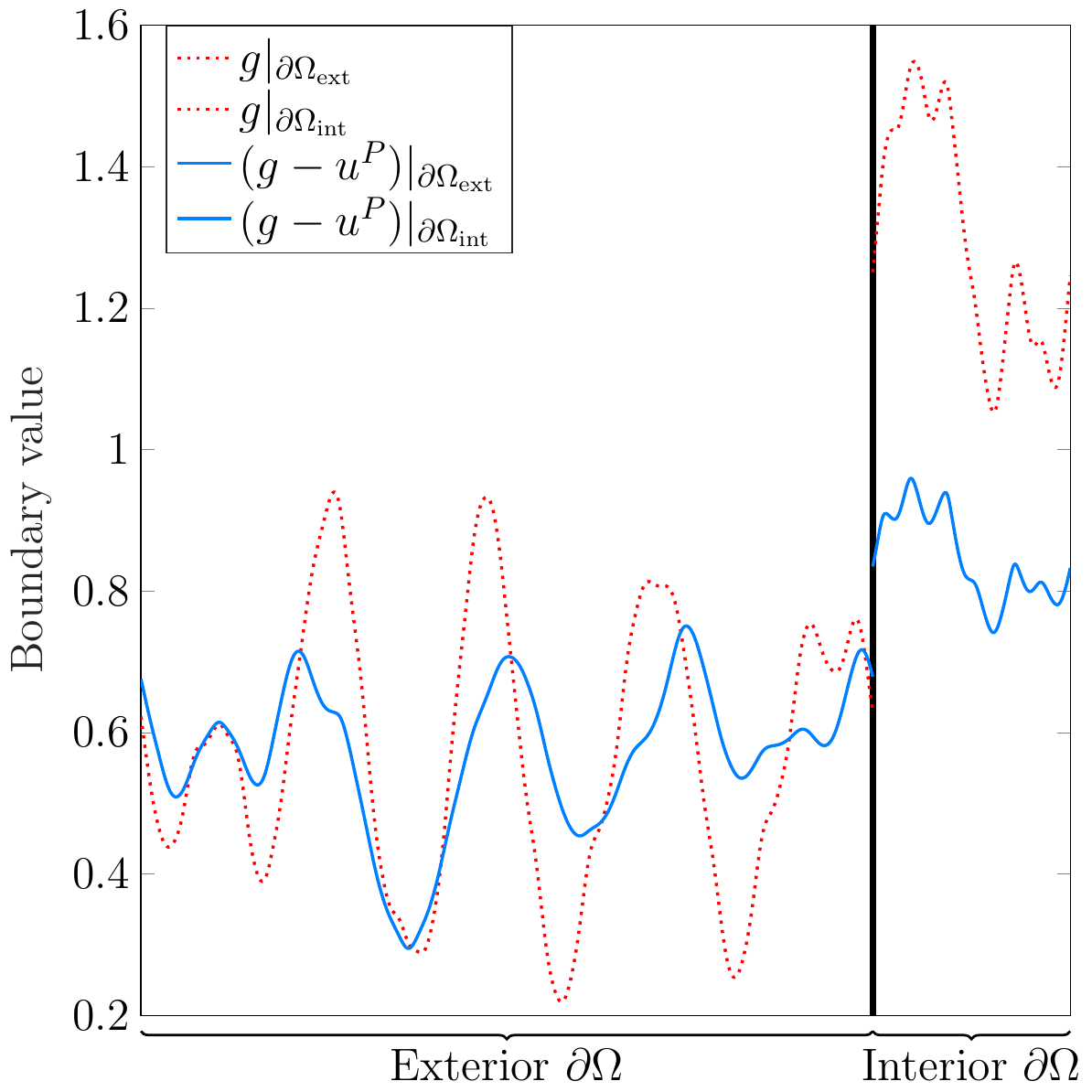}
    \end{center}
    \caption{Results for solving the Poisson equation with right hand side given by\eqref{eq:results3}. Left: Relative error as function of uniform grid resolution $N_{u}$, with $\tilde{k}$ and $M$ set by \eqref{eq:ktilde} and \eqref{eq:autoM}. Right: The boundary values $g$ and modified boundary values $g-u^{P}$ on the exterior and interior boundary segments. }
    \label{fig:function3convergenceAndbc}
\end{figure}

The numerical solution and pointwise relative error are plotted in Figure \ref{fig:example3uanderror}. The right--hand side of Figure  \ref{fig:function3convergenceAndbc} features a plot of
the error \eqref{eq:relativeError} as a function of $N_{u}$. We observe the same trend for the convergence as in the two previous examples, but require a finer grid in this setting to reach an error of  $\mathcal{O}(10^{-14})$. This is due to \eqref{eq:results3} being harder to resolve and smaller partitions are required to obtain good approximations of the interpolants. Consequently $f^{e}$ has less distance over which it goes to zero, hence a finer resolution is required to resolve it.\\

\FloatBarrier

\subsection{Performance and efficiency}
To give an idea of the computational cost of function extension by PUX, some timings are provided in Table \ref{tab:timings}, where each time is the mean of ten runs for the numerical setting presented in subsection \ref{ss:parameterSelection}. Thus the number of partitions is constant, but the grid is refined, with $\tilde{k}$ and $M$ set thereafter. For reference we have included solving the free--space Poisson equation with the method presented in subsection \ref{ss:particularSolution}. We timed:

\begin{itemize}
\item Building $A$. This is done once and the same $A$ is reused for all extension partitions. Note that RBF--QR is used whether needed or not.
\item The construction of the local extensions $\{f^{e}_{i}\}_{i=1}^{N_{p}}$. This includes: separating $A$ into $A_{\Omega}$ and $A_{E}$ as in \eqref{eq:interpol_LS_sys}, solving the least--squares problem $A_{\Omega}F=\tilde{F}_{\Omega}$ and evaluating $A_{E}F$ to obtain $\tilde{F}_{E}$ in \eqref{eq:extrapolatingf}.
\item Numerically evaluate $f^{e}$ as in \eqref{eq:functionextension}. For the extension partitions the precomputed values for $w$ are reused. For the  zero--partitions with radius other than $R_{p}$  new evaluations are required.
\item  Solve for $u^{P}$. All steps are required steps are included in the timings, even the precomputation of $\tilde{K}$, see \eqref{eq:truncatedSpectralGreen}. When excluding precomputation the time is about one tenth the tabulated result.
\end{itemize}
When timing the processes above we excluded: identifying points on the uniform grid within $R_{p}$ of the partition centre, identifying points as in $\Omega$ or $E$ and precomputing the weights $w$.\\

From Table \ref{tab:timings} we read that being able to compute $A$ once and reuse is vital for the efficiency of PUX. Solving the least--squares problem $A_{\Omega }F =\tilde{F}_{\Omega}$, see equation \eqref{eq:interpol_LS_sys}, on each partition is the most time consuming process of constructing $\{f_{i}^{e}\}_{i=1}^{N_{p}}$. This is done with a QR decomposition of $A_{\Omega}$ and backward substitution, thus scaling as $2N_{u,i}M^{2}-2/3 M^{3}+N_{u,i}M$, where $N_{u,i}$ is the number of data points in partition $\Omega_{i}$. However, for $N_{u}$ less than $400$ the local systems are not large enough for this scaling to be dominating. With downsampling to reduce $\beta_{\min}$, i.e. $N_{u,i}$, the cost of constructing $\{f^{e}_{i}\}_{i=1}^{N_{p}}$ is comparable to the cost of solving for $u^{P}$ using an efficient FFT--based solver. Observe that in our implementation $\{f^{e}_{i}\}_{i=1}^{N_{p}}$ and the weights are not computed in parallel, but doing so is trivial. The cost of the RBF--QR algorithm increases with $\varepsilon$ and the spacing between centres, both of which are comparatively large in this example. Thus these timings represent a "worst case" scenario.  
\\

\begin{table}
\centering
\caption{Timings in $[s]$ of steps for PUX and solving the free--space Poisson equation for example \ref{ss:parameterSelection}. Measurements are based on the mean of ten runs.}
\begin{threeparttable}
{\def\arraystretch{1.3}
\begin{tabular}{cccc}
\toprule
Task & $N_{u}=100$ & $N_{u}=200$ & $N_{u}=400$ \\
\midrule
Build $A$ & $0.0717$ & $0.1132$ & $0.2852$\\
Construct $\{f^{e}_{i}\}_{i=1}^{N_{p}}$ & $0.0078$ & $0.0437$ & $0.4438$\tnote{1}\\
Evaluate $f^{e}$ & $0.0043$ & $0.0090$ & $0.0307$\\
Solve for $u^{P}$ & $0.0169$ & $0.0611$ & $0.2297$\\
\bottomrule
\end{tabular}
}
\begin{tablenotes}
    \item[1] \small{Value for $\beta_{\min} \approx 19.2 $. With downsampling to $\beta_{\min} \approx 4.8$, the corresponding time is $0.191$.}
  \end{tablenotes}
\end{threeparttable}
\label{tab:timings}
\end{table}
\FloatBarrier

\section{Conclusions}
\label{s:conclusions}
In this paper, we have introduced the novel method Partition of Unity
Extension, or PUX, for numerically extending a function $f$ outside of
the complex multiply connected domain it is given on. The main
strength of the PUX method is that global regularity is obtained by
solving local problems. This is achieved by blending local
extensions on circular patches with a partition of unity 
function.  Its regularity can be chosen and
determines the global regularity. Moreover, by introducing what we
call zero partitions when defining this function, a compactly
supported function extension is obtained.\\

The performance of PUX has been thoroughly investigated by solving
the Poisson equation on multiply connected complex domains. The
particular solution comes from solving the free--space Poisson equation
for the extended $f$ given by PUX. Thereafter the Laplace equation is
solved on the given domain using a boundary integral method, with
boundary conditions modified according to the particular solution. The
final solution is the sum of the two. We have demonstrated how the various parameters are related and how to
set them, thereby significantly reducing the parameter space. By using
these guidelines the error in the solution to the Poisson equation
converges to $10^{-14}$ with an order of about $O(N_u^{-10})$, where
$N_u \times N_u$ is the total number of uniform grid points. This shows that the method described in the paper, without any additional numerical treatment, can provide precision down to round off. \\

The PUX method is simple to implement, with the exception of the
RBF--QR method used to compute the matrix $A$ as defined in section \ref{s:functionextesion}.
For this, an open--source implementation is available online
\citep{RBFQRlink}. The RBF--QR algorithm is computationally costly as compared to
the least squares solve that follows. By centring all circular
patches on a uniform grid point, the matrix $A$ however needs to be
computed only once, which yields a great reduction in the total computational
cost.  The local least squares problems  on each
patch can be solved accurately with standard methods and are typically quite small, with $100$-$400$ unknowns.\\

The PUX algorithm in this paper has been designed for when the data of
$f$ is defined on a uniform grid. Generally a finer grid is required to resolve the extension as it goes to zero, compared to resolve $f$ on $\Omega$. Since the grid is uniform the resolution is often set by the extension. It would be possible to define the
PUX method also for an adaptive grid. The implementation can be
changed to allow for patches of different sizes since this is no
restriction for the method.\\

The local extension with radial basis functions as well as the
partition of unity blending naturally extends to three
dimensions. Circular patches will become spherical patches, but
nothing conceptually changes. A paper describing the implementation of
PUX in three dimensions is forthcoming.

\section{Acknowledgements}
\label{s:acknowledgements}
This work has been supported by the Swedish Research Council under Grant No. $2015$--$04998$ and by the G\"{o}ran Gustafsson Foundation for Research in Natural Sciences and Medicine and is gratefully acknowledged.






\bibliographystyle{model1-num-names}
\bibliography{sections/func_ext}







\end{document}